\documentclass[12pt]{article}

\usepackage[english]{babel}
\usepackage[T1]{fontenc}
\usepackage{amsmath}
\usepackage{amsfonts}
\usepackage{amssymb}
\usepackage{amsthm}
\usepackage[utf8]{inputenc}
\usepackage{bm,graphicx,color,url}
\usepackage[margin=2.4cm]{geometry}
\usepackage{booktabs}
\usepackage{subfigure}
\usepackage{cite}
\usepackage{caption}
\usepackage{tikz}
\def\red{\color[rgb]{1,0,0}}
\def\magenta{\color[rgb]{1,0,1}}
\def\blue{\color[rgb]{0,0,1}}

\def\RR{\mathbb{R}}
\def\DD{{\cal D}}
\def\dd{\mathrm{d}}
\def\pmatrix{\left(\begin{array}}
\def\endpmatrix{\end{array}\right)}
\def\bfb{\bm{b}}
\def\bfc{\bm{c}}
\def\bfs{\bm{s}}
\def\bfe{\bm{e}}
\def\bfx{\bm{x}}
\def\bfeta{\bm{\eta}}
\def\bffi{\bm{\phi}}
\def\bfzero{\bm{0}}
\def\bfgamma{\bm{\gamma}}
\def\bfz{\bm{z}}
\def\bfq{\bm{q}}
\def\bfp{\bm{p}}
\def\bfw{\bm{\omega}}
\def\P{{\cal P}}
\def\I{{\cal I}}
\def\H{{\cal H}}
\def\M{{\cal M}}
\def\J{{\cal J}}
\def\U{{\cal U}}

\newtheorem{Definition}{Definition}
\newtheorem{Theorem}{Theorem}

\begin{document}
\title{Line Integral solution of Hamiltonian PDEs}

\author{Luigi Brugnano\footnote{ Dipartimento di Matematica e Informatica ``U.\,Dini'', Universit\`a di Firenze, Viale Morgagni 67/A, 
50134 Firenze, Italy; luigi.brugnano@unifi.it},\quad Gianluca Frasca-Caccia\footnote{School of Mathematics, Statistics \& Actuarial Science, University of Kent, Sibson Building,  Parkwood Road, Canterbury, CT2 7FS, UK; G.Frasca-Caccia@kent.ac.uk},\quad Felice Iavernaro\footnote{Dipartimento di Matematica, Universit\`a di Bari, Via Orabona 4, 70125 Bari, Italy; felice.iavernaro@uniba.it}\,\\[.5cm]
%\small
%$^{1}$  Dipartimento di Matematica e Informatica ``U.\,Dini'', Universit\`a di Firenze,\\ 
%\small Viale Morgagni 67/A, 
%50134 Firenze, Italy; luigi.brugnano@unifi.it\\
%\small
%$^{2}$  School of Mathematics, Statistics \& Actuarial Science, University of Kent,\\ 
%\small  Parkwood Road, Canterbury, CT2 7FS, UK; G.Frasca-Caccia@kent.ac.uk\\
%\small
%$^{3}$  Dipartimento di Matematica, Universit\`a di Bari,\\
%\small Via Orabona 4, 70125 Bari, Italy; felice.iavernaro@uniba.it
}

\maketitle

%% Author Orchid ID: enter ID or remove command
%\newcommand{\orcidauthorA}{0000-0002-6290-4107} % Add \orcidA{} behind the author's name
%\newcommand{\orcidauthorB}{0000-0002-4703-1424} % Add \orcidB{} behind the author's name
%\newcommand{\orcidauthorC}{0000-0002-9716-7370} % Add \orcidC{} behind the author's name
%
%% Authors, for the paper (add full first names)
%\Author{Luigi Brugnano $^{1}$*\orcidA{}, Gianluca Frasca-Caccia $^{2}$\orcidB{} and Felice Iavernaro $^{3}$\orcidC{}}
%
%% Authors, for metadata in PDF
%\AuthorNames{Luigi Brugnano, Gianluca Frasca-Caccia and Felice Iavernaro}
%
%% Affiliations / Addresses (Add [1] after \address if there is only one affiliation.)
%\address{%
%$^{1}$ \quad Dipartimento di Matematica e Informatica ``U.\,Dini'', Universit\`a di Firenze, Viale Morgagni 67/A, 
%50134 Firenze, Italy; luigi.brugnano@unifi.it\\
%$^{2}$ \quad School of Mathematics, Statistics \& Actuarial Science, University of Kent,   Sibson Building, Parkwood Road, Canterbury, CT2 7FS, UK; G.Frasca-Caccia@kent.ac.uk\\
%$^{3}$ \quad Dipartimento di Matematica, Universit\`a di Bari, Via Orabona 4, 70125 Bari, Italy; felice.iavernaro@uniba.it}

% Contact information of the corresponding author
%\corres{Correspondence: luigi.brugnano@unifi.it}

% Current address and/or shared authorship
%\firstnote{Current address: Affiliation 3} 
%\secondnote{All the authors contributed equally to this work.}
% The commands \thirdnote{} till \eighthnote{} are available for further notes

% Simple summary
%\simplesumm{}

% Abstract (Do not insert blank lines, i.e. \\) 
\abstract{In this paper, we report about recent findings in the numerical solution of Hamiltonian Partial Differential Equations (PDEs), by using energy-conserving line integral methods in the Hamiltonian Boundary Value Methods (HBVMs) class. In particular, we consider the semilinear wave equation, the nonlinear Schr\"odinger equation, and the Korteweg--de Vries equation, to illustrate the main features of this novel approach.\\}

% Keywords
\textbf{Keywords:} {Hamiltonian problems; energy-conserving methods; Hamiltonian Boundary Value Methods; HBVMs; line integral methods; spectral methods; Hamiltonian PDEs; semilinear wave equation; nonlinear Schr\"odinger equation; Korteweg--de Vries equation.\\}

% The fields PACS, MSC, and JEL may be left empty or commented out if not applicable
%\PACS{J0101}
 \textbf{MSC:} {65P10, 65M70, 65M20, 65L05, 65L06}
%\JEL{}
%%%%%%%%%%%%%%%%%%%%%%%%%%%%%%%%%%%%%%%%%%

%%%%%%%%%%%%%%%%%%%%%%%%%%%%%%%%%%%%%%%%%%

\section{Introduction}

The numerical solution of ordinary differential equations (ODE) problems, though researched for over sixty years, is still a very active field of investigation, following a number of trends, such as:
\begin{enumerate}
\item[a)] the search for methods suited for specific relevant classes of problems;

\item[b)] their efficient implementation on a computer;

\item[c)] the extension of existing methods to cope with wider classes of problems.
\end{enumerate}

\medskip
Point a) is particularly interesting, since it is nowadays well understood that relevant classes of problems do possess specific {\em geometric properties} in their solutions and, often, one is interested in reproducing such properties in the discrete solution obtained by a  numerical method. As matter of fact, the term {\em Geometric Integration} has been coined to denote the study of numerical methods able to preserve such properties. These latter methods, in turn, are named {\em geometric integrators}. As an example, when dealing with dissipative problems, $A$-stable methods are geometric integrators, since they retain the asymptotic stability of equilibria. Nevertheless, when stability results by first approximation do not apply, things become much more involved.  This is the case, for example, of Hamiltonian problems, i.e., problems in the form%\marginpar{\blue Hpro}
\begin{equation}\label{Hpro}
\dot y = J\nabla H(y) =: f(y), \qquad y(0) = y_0\in \mathbb{R}^{2m},  
\end{equation}
with $J=-J^\top$ and $H$ a scalar function (which we shall hereafter assume to be suitably regular), called the {\em Hamiltonian} or {\em energy}.  Due to the skew-symmetry of $J$, this latter function turns out to be conserved along the solution of (\ref{Hpro}). In fact, one has:
$$\frac{\dd}{\dd t} H(y) = \nabla H(y)^\top \dot y = \nabla H(y)^\top J\nabla H(y)=0.$$
Hamiltonian problems are very important in the applications and, for this reason, their numerical simulation has been the subject of many researches: we refer the reader, e.g., to the monographs \cite{SSC1994,LeRe2004,GNI2006,LIMbook2016,BC2016} and references therein. In particular, numerical methods able to conserve $H$ are geometric integrators, referred to as {\em energy-conserving methods}.

\medskip
Point b) is also paramount: in fact, no numerical method can be really useful, if it cannot be efficiently implemented on a computer. Therefore, a particular care has to be devoted to devise robust implementation techniques, in order to make the studied methods suitable for solving  a wide class of problems. In particular, the availability of efficient Newton-type procedures for solving the discrete problems generated by the methods turns out to be central, when numerically solving the Hamiltonian problems described at the next point.

\medskip
At last, point c) is one of the main focuses of the present paper. In fact, according to \cite[page\,157]{SSC1994}, one effective way of solving {\em Hamiltonian PDEs} is to discretize, at first, the space variable(s). In so doing, under appropriate space discretizations, one obtains a large-size  Hamiltonian problem, which can be then solved by using a suitable geometric integrator. In particular, for sake of simplicity and brevity, in this paper we shall deal with initial-boundary value problems in one space dimension, equipped with periodic boundary conditions, even though the arguments can be extended to cope with higher space dimensions, as is sketched in Section~\ref{hisp}. As was anticipated above, the numerical solution of the Hamiltonian problems arising from the space discretization of Hamiltonian PDEs will require the use of effective Newton-type procedures, in order to avoid severe step-size limitations.  

\medskip
With these premises, the present paper is devoted to report about recent findings in the numerical solution of Hamiltonian PDEs by using {\em Hamiltonian Boundary Value Methods (HBVMs)}, a class of energy-conserving Runge-Kutta methods for Hamiltonian problems. 
The novelty in their use stems from the fact that they provide effective and arbitrarily high-order energy-conserving methods for the time integration of the Hamiltonian semi-discrete problems obtained from Hamiltonian PDEs. In fact, low order methods have been mainly considered for this purpose, so far (see, e.g., \cite{CGMMOOQ2012,KS2013,GCW2014,MQ2014,GW2016,JSLW2017}). Further approaches can be found in \cite{FSS1997,DSS1998,F1999,MF2001,DLM2003,DO2011,FM2011,GX2015,FCH2018}. In more details, the structure of this paper is as follows:
\begin{itemize}
\item in Section~\ref{HBVMsec} we recall the main facts about HBVMs, also sketching their efficient {\em blended} implementation;

\item in Section~\ref{wavesec} we describe the space discretization of the semilinear wave equation,  and the efficient solution of the resulting Hamiltonian ODE problem via HBVMs. For this equation we shall provide full details, whereas the whole procedure will be only sketched for the subsequent equations;

\item in Section~\ref{NLSEsec} we see that the same approach can be used for the nonlinear Schr\"odinger equation;

\item in Section~\ref{KdVsec} we consider, instead, the Korteweg--de Vries equation;

\item Section~\ref{numsec} contains some numerical tests, aimed at showing the effectiveness of the proposed approach;

\item at last, a few conclusions are given in Section~\ref{endsec}.

\end{itemize}

\section{Hamiltonian Boundary Value Methods (HBVMs)}\label{HBVMsec}
HBVMs are energy-conserving methods derived within the framework of {\em (discrete) line integral methods}, initially proposed in \cite{IaTr2005,IaTr2006,IaPa2007,IaPa2008,IaTr2009}, and later refined in \cite{BIT2009,BIT2010,BIS2010,BIT2012_2,BIT2012_3,BIT2015}. The approach has also been extended along several directions \cite{BIT2010,BIT2010_1,BCMR2012,BI2012,BIT2012,BIT2012_1,BrSu2014,BGI2018}, including Hamiltonian BVPs \cite{ABI2015}, constrained Hamiltonian problems \cite{BGIW2018}, highly-oscillatory problems \cite{BMR2018,BIMR2018,ABI2018}, and Hamiltonian PDEs \cite{BFCI2015,BBFCI2018,BGZ2018,BIMR2018,BZL2018,BGS2019}. We also refer to the review paper \cite{BI2018} and to the monograph \cite{LIMbook2016}.

The basic idea {\em line integral methods} rely on is that the conservation of an invariant can be recast as the vanishing of a corresponding line-integral. In the case of the Hamiltonian $H$ for (\ref{Hpro}), one has:
$$H(y(t))-H(y_0) = \int_0^t \nabla H(y(\tau))^\top \dot y(\tau)\dd \tau = \int_0^t \nabla H(y(\tau))^\top J\nabla H(y(\tau))\dd \tau =0,$$
due to the fact that the integrand is identically zero. Consequently, $H(y(t))=H(y_0)$, for all $t\ge0$. Nevertheless, when dealing with a {\em discrete time} dynamics, ruled by a time-step $h>0$, one can consider a path $\sigma:[0,h]\rightarrow\RR^{2m}$ such that
\begin{equation}\label{lim1}
\sigma(0)=y_0, \qquad \sigma(h)=:y_1, \qquad y_1\approx y(h),
\end{equation}
and%\marginpar{\blue lim1,lim2}
\begin{align}\nonumber
H(y_1)-H(y_0) &= H(\sigma(h))-H(\sigma(0)) \\\label{lim2}
 &=\int_0^h \nabla H(\sigma(t))^\top \dot \sigma(t)\dd t = h\int_0^1 \nabla H(\sigma(ch))^\top \dot \sigma(ch)\dd c=0,
\end{align}
but without requiring the integrand to be identically zero. In such a case, there are infinitely many paths satisfying (\ref{lim1})--(\ref{lim2}), each providing a corresponding line integral method. In particular, we here consider a polynomial path, which we expand along the orthonormal Legendre basis:%\marginpar{\blue Leg}
\begin{equation}\label{Leg} 
P_i\in\Pi_i, \qquad \int_0^1 P_i(c)P_j(c)\dd c=\delta_{ij}. \qquad \forall i,j=0,1,\dots,
\end{equation}
where, as is usual, $\Pi_i$ is the set of polynomials of degree $i$ and $\delta_{ij}$ is the Kronecker symbol. 
In order to obtain a path $\sigma\in\Pi_s$ satisfying (\ref{lim1})-(\ref{lim2}), let us then consider the expansion%\marginpar{\blue sig1}
\begin{equation}\label{sig1}
\dot\sigma(ch) = \sum_{j=0}^{s-1}P_j(c)\gamma_j(\sigma), \qquad c\in[0,1],
\end{equation}
in terms of the $s$ unknown vector coefficients $\{\gamma_j(\sigma)\}$. In order to fulfill (\ref{lim1}), integrating both sides of (\ref{sig1}) and taking into account that (see (\ref{Leg})) $\int_0^1P_j(c)\dd c=\delta_{j0}$, one obtains:%\marginpar{\blue sig}
\begin{equation}\label{sig}
\sigma(ch) = y_0+h\sum_{j=0}^{s-1} \int_0^c P_j(\tau)\dd \tau\gamma_j(\sigma), \quad c\in[0,1], \quad \Rightarrow\quad y_1\equiv \sigma(h)=y_0+h\gamma_0(\sigma).
\end{equation}
Taking into account (\ref{sig1}), condition (\ref{lim2}) becomes
\begin{eqnarray*}
\int_0^1 \nabla H(\sigma(ch))^\top \dot \sigma(ch)\dd c&=&\int_0^1 \nabla H(\sigma(ch))^\top \sum_{j=0}^{s-1} P_j(c)\gamma_j(\sigma)\dd c\\
&=&\sum_{j=0}^{s-1}\left( \int_0^1P_j(c)\nabla H(\sigma(ch))\dd c\right)^\top\gamma_j(\sigma) ~=~0,
\end{eqnarray*}
which is satisfied by choosing (see (\ref{Hpro})): 
\begin{equation}\label{gammaj}
\gamma_j(\sigma) = J\int_0^1 P_j(c)\nabla H(\sigma(ch))\dd c \equiv \int_0^1 P_j(c) f(\sigma(ch))\dd c,
\end{equation}
because of the skew-symmetry of matrix $J$. Therefore, this specific energy-conserving line integral method is defined by the polynomial path $\sigma$, whose coefficients satisfy the following set of $s$ nonlinear vector equations, derived from (\ref{sig}) and (\ref{gammaj}): 
\begin{equation}\label{dispro}
\gamma_j(\sigma) = \int_0^1 P_j(c) f\left(y_0+h\sum_{i=0}^{s-1} \int_0^c P_i(\tau)\dd \tau\gamma_i(\sigma)\right)\dd c, \qquad j=0,\dots,s-1.
\end{equation}
Moreover, it can be proved that  $\sigma(h)-y(h)=O(h^{2s+1})$, i.e., the approximation procedure has order $2s$ \cite[Theorem\,1]{BIT2012_3} (see also \cite{BI2018}). However, this procedure does not yet provide a numerical method since, quoting e.g. Dahlquist and Bj\"ork \cite[page\,521]{DaBj2008}, ``{\em as is well known, even many relatively simple integrals cannot be expressed in finite terms of elementary functions, and thus must be evaluated by numerical methods}.'' In particular, since we are dealing with a polynomial approximation, we consider the Gaussian interpolatory quadrature rule, based at the zeros $0<c_1<\dots<c_k<1$ of $P_k$, whose weights we denote, respectively, by $b_1,\dots,b_k$, which is well-known to have order $2k$.\footnote{ I.e., it is exact for polynomial integrands up to order $2k-1$.} Consequently, with reference to (\ref{gammaj}), we obtain the approximation 
\begin{equation}\label{hgammaj}
\gamma_j(\sigma) \approx \sum_{\ell=1}^k b_\ell P_j(c_\ell) f(\sigma(c_\ell h)) =: \hat\gamma_j, \qquad j=0,\dots,s-1,
\end{equation}
where, for sake of brevity, we continue to denote $\sigma$ the polynomial approximation. The new discrete problem is then given by
\begin{equation}\label{dispro1}
\hat\gamma_j = \sum_{\ell=1}^k b_\ell P_j(c_\ell) f\left(y_0+h\sum_{i=0}^{s-1} \int_0^{c_\ell} P_i(\tau)\dd \tau\hat\gamma_i\right), \qquad j=0,\dots,s-1,
\end{equation} 
which remarkably has, alike (\ref{dispro}), dimension $s$, independently of $k$. 
\begin{Definition}\label{hbvmdef}
 The discrete problem (\ref{dispro1}) defines a HBVM$(k,s)$ method. The limit as $k\rightarrow\infty$, given by (\ref{dispro}), defines a HBVM$(\infty,s)$ formula.
\end{Definition}

It is possible to prove the following result \cite{LIMbook2016,BI2018}. 

\begin{Theorem}\label{hbvmth}
For all $k\ge s$, by using the $k$ Gauss-Legendre abscissae, a HBVM$(k,s)$ method is symmetric and of order $2s$. Moreover, it reduces to the $s$-stage Gauss collocation method, when $k=s$. Concerning energy-conservation when applied for solving (\ref{Hpro}),  one has: 
\begin{equation}\label{Hcons}
H(y_1)-H(y_0) = \left\{ \begin{array}{cc} 0, &\mbox{if}\quad H\in\Pi_\nu \quad\mbox{with}\quad \nu\le 2k/s,\\[2mm]
O(h^{2k+1}), &\mbox{otherwise}.\end{array}\right.
\end{equation}
\end{Theorem}

It is worth mentioning that, because of (\ref{Hcons}), by choosing $k$ large enough one can either obtain:
\begin{itemize}
\item an {\em exact} conservation of energy, when $H$ is a polynomial;
\item a {\em practical} conservation of energy, otherwise. In fact, in such a case, it is enough that the energy error falls within the round-off error level. 
\end{itemize}

\subsection{Runge-Kutta form of HBVM$(k,s)$}
It is possible to see that, actually, a HBVM$(k,s)$ method is a $k$-stage Runge-Kutta method. In fact, by setting in (\ref{hgammaj}) $Y_\ell:=\sigma(c_\ell h)$, $\ell=1,\dots,k$, one obtains:%\marginpar{\blue rk1}
\begin{eqnarray}\nonumber
Y_i &\equiv& \sigma(c_ih) ~=~ y_0 + h\sum_{j=0}^{s-1}\int_0^{c_i} P_j(\tau)\dd\tau\hat\gamma_j ~=~y_0 + h\sum_{j=0}^{s-1}\int_0^{c_i} P_j(\tau)\dd\tau\sum_{\ell=1}^k b_\ell P_j(c_\ell) f(Y_\ell)\\ \label{rk1}
&=& y_0+h \sum_{j=1}^k \left[ b_j \sum_{\ell=0}^{s-1} \int_0^{c_i} P_\ell(\tau)\dd\tau \,P_\ell(c_j)\right] f(Y_j), \qquad i=1,\dots,k,
\end{eqnarray}
with the new approximation given by%\marginpar{\blue rk2}
\begin{equation}\label{rk2}
y_1 = y_0 + h\hat\gamma_0 \equiv y_0 + h\sum_{i=1}^k b_i f(Y_i).
\end{equation}
It can be readily seen that (\ref{rk1})-(\ref{rk2}) define the $k$-stage Runge-Kutta method with Butcher tableau%\marginpar{\blue Butab}
\begin{equation} \label{Butab}
\begin{array}{c|c} \bfc & \I_s \P_s^\top\Omega \\ \hline & \bfb^\top\end{array} 
\end{equation}
with%\marginpar{\blue cbOm}
\begin{equation}\label{cbOm}
\bfc = \pmatrix{c} c_1 \\ \vdots \\ c_k\endpmatrix, \qquad
\bfb = \pmatrix{c} b_1 \\ \vdots \\ b_k\endpmatrix, \qquad
\Omega = \pmatrix{ccc} b_1 \\ &\ddots \\ && b_k\endpmatrix, 
\end{equation}
and%\marginpar{\blue PI}
\begin{equation}\label{PI}
\P_s = \pmatrix{ccc} P_0(c_1) & \dots & P_{s-1}(c_1)\\
\vdots & &\vdots\\
P_0(c_k) & \dots & P_{s-1}(c_k)
\endpmatrix, ~ \I_s = \pmatrix{ccc} \int_0^{c_1}P_0(x)\dd x & \dots & \int_0^{c_1} P_{s-1}(x)\dd x\\
\vdots & &\vdots\\
\int_0^{c_k}P_0(x)\dd x & \dots & \int_0^{c_k} P_{s-1}(x)\dd x
\endpmatrix ~\in\RR^{k\times s}.
\end{equation}
For this Runge-Kutta method, the stage equation (\ref{rk1}) has (block) dimension $k$ and is given by
$$Y = \bfe\otimes y_0 + h\I_s\P_s^\top\Omega\otimes I_{2m} f(Y), \quad Y = \pmatrix{c} Y_1\\ \vdots\\ Y_k\endpmatrix, \,
 f(Y) = \pmatrix{c} f(Y_1)\\ \vdots\\ f(Y_k)\endpmatrix, \, \bfe = \pmatrix{c} 1\\ \vdots\\ 1\endpmatrix\in\RR^k,$$
 having set, in general, $I_r\in\RR^{r\times r}$  the identity matrix. 
Nonetheless, the equivalent discrete problem (\ref{dispro1}), whose dimension is $s$ independently of $k$,
turns out to be given by:
\begin{equation}\label{dispro2}
F(\hat\bfgamma) := \hat\bfgamma - \P_s^\top\Omega\otimes I_{2m} f\left( \bfe\otimes y_0 + h\I_s\otimes I_{2m}\hat\bfgamma\right) = \bfzero,
\end{equation} where
\begin{equation}\label{hgamma}
\hat\bfgamma = \pmatrix{c} \hat\gamma_0\\ \vdots \\ \hat\gamma_{s-1}\endpmatrix, \qquad \hat\gamma_i\in\RR^{2m},\quad i=0,\dots,s-1. 
\end{equation}
Once (\ref{dispro2}) is solved, according to (\ref{rk2}) the new approximation is given by~ $y_1=y_0+h\hat\gamma_0$.

\subsection{Special second-order problems}
Sometimes, the problem (\ref{Hpro}) assumes the form of a {\em special second-order problem}, 
\begin{equation}\label{Hpro2}
\ddot q = \nabla U(q), \qquad q(0)=q_0,~\dot q(0)=p_0\in\RR^m,
\end{equation}
for which, setting $p=\dot q$,\, $y=\pmatrix{c} q\\ p\endpmatrix$\, and\, $H(y)\equiv H(q,p) = \frac{1}2 p^\top p-U(q)$. In such a case, the dimension of the blocks of the discrete problem can be halved. In fact, by using (\ref{Butab}) for solving (\ref{Hpro2}),  one sees that the stage equations for $q$ and $p$ are respectively given by: 
\begin{equation}\label{QP}
Q = \bfe\otimes q_0 +h\I_s\P_s^\top\Omega \otimes I_m P, \qquad P = \bfe\otimes p_0 +h\I_s\P_s^\top\Omega \otimes I_m \nabla U(Q),
\end{equation}
having set
$$Q = \pmatrix{c} Q_1\\ \vdots\\ Q_k\endpmatrix, \qquad P=\pmatrix{c} P_1\\ \vdots\\ P_k\endpmatrix,\qquad \nabla U(Q) = \pmatrix{c} \nabla U(Q_1)\\ \vdots\\ \nabla U(Q_k)\endpmatrix.$$
Plugging the second equation in (\ref{QP}) into the first one,  considering that
 $\I_s\P_s^\top\Omega\bfe = \bfc$ and, moreover, 
\begin{equation}\label{Xs}
\P_s^\top\Omega\I_s = X_s \equiv \pmatrix{cccc} 
\xi_0 & -\xi_1\\
\xi_1 &0 &\ddots\\
         &\ddots &\ddots &-\xi_{s-1}\\
         &           &\xi_{s-1}&0\endpmatrix, \quad \xi_i = \left(2\sqrt{|4i^2-1|}\right)^{-1}, \quad i=0,\dots,s-1,
         \end{equation}
one then obtains:%\marginpar{\blue Q}
\begin{equation}\label{Q}
Q = \bfe\otimes q_0 +h\bfc\otimes p_0 + h^2\I_s X_s \P_s^\top\Omega\otimes I_m\nabla U(Q).
\end{equation}
Setting (compare with (\ref{hgamma}))%\marginpar{\blue bgamma}  
\begin{equation}\label{bgamma}
\bar\bfgamma \equiv \pmatrix{c}\bar\gamma_0\\ \vdots \\ \bar\gamma_{s-1}\endpmatrix = \P_s^\top\Omega\otimes I_m\nabla U(Q),
\qquad \bar\gamma_i\in\RR^m,\quad i=0,\dots,s-1,
\end{equation}
and taking into account (\ref{Q}), one then obtains the new discrete problem (compare with (\ref{dispro2})):%\marginpar{\blue dispro3}
\begin{equation}\label{dispro3}
G(\bar\bfgamma) := \bar\bfgamma-\P_s^\top\Omega\otimes I_m \nabla U\left( \bfe\otimes q_0+h\bfc\otimes p_0 + h^2\I_s X_s \otimes I_m \bar\bfgamma\right)=\bfzero.
\end{equation}
Once it has been solved, it can be seen that the new approximations are given by (see, e.g., \cite[Chapter\,4]{LIMbook2016}):
$$q_1 = q_0 + h p_0 + h^2\left( \xi_0\bar\gamma_0 -\xi_1\bar\gamma_1\right), \qquad p_1 = p_0 + h\bar\gamma_0,$$
where $\xi_0$ and $-\xi_1$ are the nonzero entries on the first row of matrix $X_s$ defined in (\ref{Xs}). 

\subsection{Blended iteration} The efficient solution of the discrete problem (\ref{dispro2}) has been studied in a series of papers \cite{BIT2011,BFCI2014,BFCI2014_1,LIMbook2016}. We here recall the main facts about the so called {\em blended implementation} of HBVMs, which represents a Newton-type iteration for solving (\ref{dispro2}). This approach, at first sketched in \cite{Br2000}, has then been analyzed in \cite{BrMa2002} and developed in \cite{BrMa2007,BrMa2009,BrMa2009_1}. It has been then implemented in the Fortran codes {\tt BiM} \cite{BrMa2004} and {\tt BiMD} \cite{BrMaMu2006}, for the numerical solution of stiff ODE-IVPs and linearly implicit DAEs: both codes can be retrieved at \cite{BIM}; the latter code is also available at the {\em Test Set for IVP Solvers} \cite{testset}. The blended implementation of HBVMs has then been considered in \cite{BIT2011} and implemented in the Matlab function {\tt hbvm} available at the url \cite{LIMbook}. We also mention that, more recently, this approach has been also considered for RKN methods \cite{WMF2017}. 

Let us then consider the simplified Newton iteration for solving (\ref{dispro2}) which, by taking into account (\ref{Xs}) amounts to solving the following set of linear systems:%\marginpar{\blue blend1}
\begin{equation}\label{blend1}
\left[ I_s\otimes I_{2m} - hX_s\otimes f'(y_0)\right] \Delta \hat\bfgamma^\ell = -F(\hat\bfgamma^\ell), \qquad 
\ell=0,1,\dots,\end{equation}
with $f'(y_0)$ the Jacobian of $f$ evaluated at $y_0$. 
This iteration, though straightforward and very effective, requires, however, the factorization of a $2ms\times 2ms$ matrix, which can be cumbersome, when $s$ and/or $m$ are large. To get rid of this problem, by considering that matrix $X_s$ is nonsingular one at first considers the following equivalent formulation of (\ref{blend1}), having set $\rho_s$ a positive, and for the moment unspecified, parameter:%\marginpar{\blue blend2}
\begin{equation}\label{blend2}
\rho_s\left[ X_s^{-1}\otimes I_{2m} - hI_s\otimes f'(y_0)\right] \Delta \hat\bfgamma^\ell = -(\rho_s X_s^{-1}\otimes I_{2m})F(\hat\bfgamma^\ell), \qquad \ell=0,1,\dots.\end{equation}
The next step is to consider the {\em blending} of the two equivalent formulations (\ref{blend1}) and (\ref{blend2}) with weights $\theta_s$ and $I_s\otimes I_{2m} - \theta_s$, respectively, where:
\begin{equation}\label{tetas}
\theta_s = I_s\otimes \Sigma^{-1}, \qquad \Sigma = [I_{2m}-h\rho_s f'(y_0)].
\end{equation}
In so doing, one obtains a new linear system, whose coefficient matrix has the inverse which can be approximated by $\theta_s$. Skipping the details (for which we refer to \cite{BIT2011}, see also \cite{LIMbook2016,BI2018}), one then obtains the following {\em blended iteration} for solving (\ref{dispro2}):
\begin{equation}\label{blend3}
\bfeta^\ell=-F(\hat\bfgamma^\ell), \quad \bfeta_1^\ell = \left(\rho_sX_s^{-1}\otimes I_{2m}\right)\bfeta^\ell, \quad \Delta\hat\bfgamma^\ell = \theta_s\left[\bfeta_1^\ell+\theta_s\left(\bfeta^\ell-\bfeta_1^\ell\right)\right], \, \ell=0,1,\dots,
\end{equation}
which only requires to factor the matrix $\Sigma$ in (\ref{tetas}), having the same size as that of the continuous problem. Concerning the choice of the parameter $\rho_s$, as is shown in \cite{BrMa2002}, the optimal choice, based on a linear convergence analysis, turns out to be:
\begin{equation}\label{ros}
\rho_s = \min_{\lambda\in\sigma(X_s)} |\lambda|,
\end{equation}
where, as is usual, $\sigma(X_s)$ is the spectrum of $X_s$.

In the case of the special second-order problem (\ref{Hpro2}), the simplified Newton iteration for solving (\ref{dispro3}) becomes: 
\begin{equation}\label{blend4}
\left[ I_s\otimes I_m - h^2X_s^2\otimes \nabla^2 U(q_0)\right] \Delta \bar\bfgamma^\ell = -G(\bar\bfgamma^\ell), \qquad 
\ell=0,1,\dots,\end{equation}
with $\nabla^2 U(q_0)$ the Hessian of $U$ evaluated at $q_0$. Consequently, similar steps as above can be repeated, via the following formal substitutions:
$$F\rightarrow G, \quad \hat\bfgamma\rightarrow\bar\bfgamma, \quad f'(y_0)\rightarrow \nabla^2U(q_0), \quad I_{2m}\rightarrow I_m, \quad h\rightarrow h^2, \quad X_s\rightarrow X_s^2, \quad \rho_s\rightarrow\rho_s^2.$$ As a result, the blended iteration for solving (\ref{dispro3}) is given by:
\begin{equation}\label{blend5}
\bfeta^\ell=-G(\bar\bfgamma^\ell), \quad \bfeta_1^\ell = \left(\rho_s^2X_s^{-2}\otimes I_m\right)\bfeta^\ell, \quad \Delta\bar\bfgamma^\ell = \theta_s\left[\bfeta_1^\ell+\theta_s\left(\bfeta^\ell-\bfeta_1^\ell\right)\right], \quad \ell=0,1,\dots,
\end{equation}
with the parameter $\rho_s$ still given by (\ref{ros}) and 
\begin{equation}\label{tetas1}
\theta_s = I_s\otimes \Sigma^{-1}, \qquad \Sigma = [I_m-h^2\rho_s^2 \nabla^2 U(q_0)].
\end{equation}
Consequently, also in such a case, one has only to factor a matrix having the same size as that of the continuous problem.

\subsection{Blended iteration for semilinear problems}
Once more, we stress that the availability of a Newton-type iteration for solving (\ref{dispro2}) is paramount, in order to avoid severe step-size limitations, when such a problem is derived from the space discretization of Hamiltonian PDEs. In fact, in such a case, the resulting ODE problem turns out to be in the form
\begin{equation}\label{semilin}
\dot y = Ay + g(y), \qquad y(0)=y_0\in\RR^{2m},
\end{equation}
with the dimension and the norm of matrix $A$ tending to infinity, as the space discretization is made more and more accurate, whereas $\|g\|$ remains bounded, if the solution is bounded. Consequently,  one can consider a constant approximation of the Jacobian of the right-hand side of (\ref{semilin}), given by the matrix $A$ of the linear term. As a result, the matrix $\Sigma$ defined in (\ref{tetas}) becomes
\begin{equation}\label{Sigma1}
\Sigma = I_{2m}-h\rho_s A,
\end{equation}
which is {\em constant for all time-steps} and, consequently, it needs to be factored only once.  

Similarly, when problem (\ref{Hpro2}) is in the form
\begin{equation}\label{semilin2}
\ddot q = -A^2q + g(q), \qquad q(0)=q_0, ~\dot q(0)=p_0\in\RR^m,
\end{equation}
with $A^2$ symmetric and semi-positive definite, and $\|A^2\|\gg \|g\|$, one can approximate the matrix $\Sigma$ in (\ref{tetas1}) as
\begin{equation}\label{Sigma2}
\Sigma = I_{m}+h^2\rho_s^2 A^2,
\end{equation}
which,  also in this case, is {\em constant for all time-steps} and needs to be factored only once.  

We end this subsection by stressing that, for the problems that we shall consider in the sequel, matrix $A$ in (\ref{Sigma1}), or matrix $A^2$ in (\ref{Sigma2}), has a block structure with {\em diagonal blocks}. As a result, the corresponding blended iterations (\ref{blend3}) and (\ref{blend5}) are computationally inexpensive. Moreover, the linear algebra can be made still more efficient, as is done in the Matlab function {\tt hbvm} available at \cite{LIMbook}, by considering a matrix formulation of the iteration \cite{LIMbook2016,Si2016}.

\subsection{HBVMs as spectral methods in time}

To conclude this quick introduction to HBVMs, we mention their use as {\em spectral methods in time}, which has been the subject of recent investigations \cite{BMR2018,BIMR2018,ABI2018}. We mention that the use of Runge-Kutta methods as spectral methods in time has been considered previously in  
\cite{B1997,BS2000,BS2000_1,TC2007} (see also \cite{BIT2012_3}). In more details, if we consider the expansion of the right-hand side of (\ref{Hpro}), on the interval $[0,h]$, along the Legendre basis (\ref{Leg}), one has: 
\begin{equation}\label{infty}
\dot y(ch) = f(y(ch)) \equiv \sum_{j\ge 0} P_j(c) \gamma_j(y), \qquad c\in[0,1],
\end{equation}
where $\gamma_j(y)$ is defined according to (\ref{gammaj}), by formally replacing $\sigma$ by $y$. On the other hand, the polynomial approximation $\sigma$ defined in (\ref{sig1}) is obtained by truncating the previous series after $s$ terms. However, by considering that  
$$\int_0^1 \|f(y(ch))\|_2^2 \dd c = \sum_{j\ge0} \|\gamma_j(y)\|_2^2,$$
one has that 
$$\|\gamma_j(y)\|_2\rightarrow 0, \qquad j\rightarrow\infty,$$
the more regular $f(y)$, the faster the convergence to 0 of $\|\gamma_j(y)\|_2$, as $j\rightarrow\infty$. Consequently, when using a finite precision arithmetic with machine epsilon $\varepsilon$, if one truncates the expansion (\ref{infty}) when the Fourier coefficient $\gamma_s(y)$ is negligible, w.r.t. the previous ones, then one obtains that (\ref{infty}) and (\ref{sig1}) become indistinguishable, in the used finite precision arithmetic. A straightforward criterion for this to happen,  considered in \cite{BMR2018,BIMR2018}, is to require that
\begin{equation}\label{tol}
\|\gamma_s(y)\|_2 < tol\cdot \max_{j=0,\dots,s-1} \|\gamma_j(y)\|_2,  
\end{equation}
with $tol\approx \varepsilon$. Moreover,  the analysis in \cite{ABI2018} shows that one could even use $tol\approx \sqrt{\varepsilon}$ in (\ref{tol}), still obtaining full machine accuracy at $t=h$. At last (see (\ref{hgammaj})), by choosing $k$ large enough, one may obtain full machine accuracy in the approximation of $\gamma_j(\sigma)$ by means of $\hat\gamma_j$, $j=0,\dots,s-1$. As a result, the use of HBVMs as spectral methods in time (which we shall denote by SHBVMs, as an abbreviation for {\em spectral HBVMs}) usually requires the use of relatively large values of $s$ and $k$. This, in turn, is not a big issue; in fact:
\begin{itemize}
\item on one hand, we have the availability of the blended iteration (\ref{blend3}) (or (\ref{blend5})), whose computational cost is mildly affected by such parameters, also considering the approximation (\ref{Sigma1}) (or (\ref{Sigma2}));

\item on the other hand, SHBVMs will allow the use of relatively large time-steps.
\end{itemize}
Summing all up, overall SHBVMs will result to be extremely effective and competitive, as is testified by the numerical tests reported in Section~\ref{numsec} (see also \cite{ABI2018,BMR2018,BIMR2018}).
 
%%%%%%%%%%%%%%%%%%%%%%%%%%%%%%%%%%%%%%%%%%
\section{The semilinear wave equation}\label{wavesec}

The first Hamiltonian PDE that we consider is the semilinear wave equation:%\marginpar{\blue wave}
\begin{eqnarray}\label{wave}
u_{tt}(x,t) &=& u_{xx}(x,t) - f'(u(x,t)), \qquad (x,t)\in[a,b]\times[0,T], \\
u(x,0)&=&u_0(x), \quad u_t(x,0)~=~v_0(x), \qquad x\in[a,b],\nonumber\end{eqnarray}
with $f'$ the derivative of $f$. The problem (\ref{wave}) is completed by prescribing periodic boundary conditions. Hereafter, we shall assume the solution to be suitably regular, as a periodic function in space. Moreover, for sake of brevity, we shall omit the arguments of the involved functions, when not necessary. By setting $v=u_t$, one obtains that (\ref{wave}) is a Hamiltonian PDE, with {\em Hamiltonian functional}%\marginpar{\blue waveH}
\begin{equation}\label{waveH}
\H[u,v](t) = \frac{1}2\int_a^b \left[ v^2(x,t)+u_x^2(x,t) +2f(u(x,t))\right]\dd x =: \int_a^b L(x,t,u,u_x,v)\dd x,
\end{equation} 
so that, by setting 
$$\nabla\H = \pmatrix{c} \delta_u \H\\ \delta_v \H\endpmatrix,$$
the vector of the {\em functional derivatives} of $\H$, with 
$$\delta_u\H = (\partial_u -\partial_x\partial_{u_x})L \equiv f'(u)-u_{xx}, \qquad \delta_v \H = \partial_v L \equiv v,$$
one has:%\marginpar{\blue J2}
\begin{equation}\label{J2}
\pmatrix{c} u_t\\ v_t\endpmatrix = J_2 \nabla\H,\qquad J_2:=\pmatrix{cc} &1\\ -1\endpmatrix,
\end{equation}
which is formally in the form (\ref{Hpro}). As in the ODE case, also now one has the conservation of the Hamiltonian.%\marginpar{\red waveth}

\begin{Theorem}\label{waveth}
Assuming that the solution of (\ref{wave}) is suitably smooth in space, the Hamitonian (\ref{waveH}) is conserved, when periodic boundary conditions are prescribed.
\end{Theorem}
\begin{proof}
In fact, from (\ref{wave}) and (\ref{waveH}), and taking into account that $v=u_t$, one has:
\begin{eqnarray*}
\dot\H[u,v] &=& \int_a^b L_t\dd x = \int_a^b \left[v v_t + u_x u_{xt}+f'(u)u_t\right]\dd x = \int_a^b \left[u_xv_x +v( v_t+f'(u))\right]\dd x \\
&=& \int_a^b \left[u_xv_x+vu_{xx}\right]\dd x  ~=~ \left[ u_xv\right]_{x=a}^{x=b} ~=~0,
\end{eqnarray*}
because of the periodic boundary conditions.
\end{proof}

In order to numerically solve (\ref{wave}), according to what sketched in the introduction, we at first discretize the space variable, with the aim of obtaining a corresponding Hamiltonian ODE problem. For this purpose, we consider the following orthonormal basis on the interval $[a,b]$, which takes into account of the periodic boundary conditions \cite{BFCI2015,LIMbook2016,BBFCI2018,BGS2019,BGZ2018,BI2018}:%\marginpar{\blue c0, cj, sj}
\begin{eqnarray}\label{c0}
c_0(x) &\equiv& (b-a)^{-\frac{1}2},\\ \label{cj}
c_j(x) &=& \sqrt{\frac{2}{b-a}}\cos\left(2\pi j\frac{x-a}{b-a}\right), \qquad x\in[a,b],\\
s_j(x) &=& \sqrt{\frac{2}{b-a}}\sin\left(2\pi j\frac{x-a}{b-a}\right), \qquad j=1,2,\dots.\label{sj}
\end{eqnarray}
In fact, for all allowed $i,j$, one has:%\marginpar{\blue ortoab}
\begin{equation}\label{ortoab}
\int_a^b c_i(x)c_j(x)\dd x =\delta_{ij} =\int_a^b s_i(x)s_j(x)\dd x, \qquad \int_a^b c_i(x)s_j(x)\dd x =0.
\end{equation} 
Consequently, for suitable time dependent coefficients $\alpha_0(t),\alpha_1(t),\beta_1(t),\dots$, one has:%\marginpar{\blue expu}
\begin{equation}\label{expu}
u(x,t) = c_0(x)\alpha_0(t)+\sum_{j\ge 1} \left[c_j(x)\alpha_j(t)+s_j(x)\beta_j(t)\right].
\end{equation}
The infinite expansion (\ref{expu}) can be cast in vector form, by defining the infinite-dimensional vectors%\marginpar{\blue wq}
\begin{equation}\label{wq}
\bfw(x) = \pmatrix{cccc} c_0(x), &s_1(x), &c_1(x), &\dots\endpmatrix^\top,\quad
\bfq(t) = \pmatrix{cccc} \alpha_0(t), &\beta_1(t), &\alpha_1(t), &\dots\endpmatrix^\top,
\end{equation}
as%\marginpar{\blue expu1}
\begin{equation}\label{expu1}
u(x,t) = \bfw(x)^\top\bfq(t).
\end{equation}
By also introducing the infinite matrix%\marginpar{\blue Dmat}
\begin{equation}\label{Dmat}
D = \left(\frac{2\pi}{b-a}\right) \pmatrix{cccc} 0\\ &1\cdot I_2 \\ &&2\cdot I_2 \\ &&&\ddots\endpmatrix,
\end{equation}
and considering that (\ref{ortoab}) can be written in matrix form as%\marginpar{\blue ortoabmat}
\begin{equation}\label{ortoabmat}
\int_a^b \bfw(x)\bfw(x)^\top \dd x = I,
\end{equation}
the identity operator, we then prove the following result.%\marginpar{\red waveth2}

\begin{Theorem}\label{waveth2} With reference to (\ref{c0})--(\ref{ortoabmat}), problem (\ref{wave}) can be rewritten as the special second-order problem%\marginpar{\blue wave1}
\begin{eqnarray}\label{wave1}
\ddot \bfq &=& -D^2\bfq -\int_a^b \bfw(x)f'(\bfw(x)^\top\bfq)\dd x, \qquad t\in[0,T], \\ \nonumber
\bfq(0) &=&\int_a^b \bfw(x)u_0(x)\dd x ~=:~ \bfq_0, \qquad \dot\bfq(0) ~=~\int_a^b \bfw(x)v_0(x)\dd x~=:~\bfp_0.
\end{eqnarray}
By setting $\bfp=\dot\bfq$, this latter problem is Hamiltonian, with Hamiltonian function%\marginpar{\blue waveH1}
\begin{equation}\label{waveH1}
H(\bfq,\bfp) = \frac{1}2\left( \bfp^\top \bfp + \bfq^\top D^2\bfq + 2\int_a^b f(\bfw(x)^\top\bfq)\dd x\right),
\end{equation}
which turns out to be equivalent to the Hamiltonian functional (\ref{waveH}).
\end{Theorem}
\begin{proof} Problem (\ref{wave1}) is clearly Hamiltonian, w.r.t. the Hamiltonian (\ref{waveH1}), since
$$\dot\bfq = \partial_{\bfp} H(\bfq,\bfp), \qquad \dot\bfp = -\partial_{\bfq} H(\bfq,\bfp).$$
Let us then show that:
\begin{itemize}
\item (\ref{wave1}) is equivalent to (\ref{wave});
\item (\ref{waveH1}) is equivalent to (\ref{waveH}).
\end{itemize}

Concerning the first point, we observe that
$$u_{tt}(x,t) = \bfw(x)^\top\ddot\bfq(t), \qquad u_{xx}(x,t) = \bfw''(x)^\top\bfq(t)\equiv -\bfw(x)^\top D^2 \bfq(t),$$
with an obvious meaning of $\bfw"(x)$, so that (\ref{wave}) can be rewritten as
$$\bfw(x)^\top\ddot\bfq(t) = -\bfw(x)^\top D^2 \bfq(t) - f'(\bfw(x)^\top\bfq(t)), \qquad (x,t)\in[a,b]\times[0,T].$$
Multiplying both sides by $\bfw(x)$, then integrating in space from $a$ to $b$, and taking into account (\ref{ortoabmat}), give us (\ref{wave1}). 

Concerning the second point, the statement easily follows by considering that $v(x,t)=u_t(x,t)=\bfw(x)^\top\bfp$, and
$$\bfp^\top\bfp = \int_a^b \bfp^\top\bfw(x)\bfw(x)^\top\bfp\dd x = \int_a^b \dot\bfq^\top\bfw(x)\bfw(x)^\top\dot\bfq\dd x = \int_a^b v^2\dd x.$$ Moreover, by defining the matrix%\marginpar{\blue bDmat}
\begin{equation}\label{bDmat}
\bar D = \left(\frac{2\pi}{b-a}\right) \pmatrix{cccc} 0\\ &1\cdot J_2 \\ &&2\cdot J_2 \\ &&&\ddots\endpmatrix = -\bar D^\top,
\end{equation}
where matrix $J_2$ is that defined in (\ref{J2}), one has:
$$u_x(x,t) = \bfw'(x)^\top\bfq(t) \equiv \left[ \bar D\bfw(x)\right]^\top\bfq(t), \qquad  \bar D\bar D^\top = D^2,$$
so that, by taking again into account (\ref{ortoabmat}), one obtains:
$$\bfq^\top D^2\bfq = \bfq^\top\bar D\bar D^\top\bfq = \int_a^b \bfq^\top\bar D \bfw(x)\bfw(x)^\top\bar D^\top\bfq\dd x = \int_a^b u_x^2\dd x.$$
The proof is completed by considering that, from (\ref{expu1}), $f(\bfw(x)^\top\bfq(t)) = f(u(x,t))$.
\end{proof}

\subsection{Discretization} In order to solve problem (\ref{wave1}) on a computer, the infinite expansion (\ref{expu}) must be truncated at a convenient index $N$. In so doing, (\ref{wq}) and (\ref{Dmat}) respectively become%\marginpar{\blue wqD1}
\begin{equation}\label{wqD1}
\bfw(x) = \pmatrix{c} c_0(x) \\ s_1(x) \\c_1(x) \\ \vdots \\ s_N(x)\\ c_N(x)\endpmatrix,\,
\bfq(t) = \pmatrix{c} \alpha_0(t) \\ \beta_1(t)  \\ \alpha_1(t) \\ \vdots \\ \beta_N(t)  \\ \alpha_N(t)\endpmatrix,\,
D = \left(\frac{2\pi}{b-a}\right) \pmatrix{ccccc} 0\\ &1\cdot I_2 \\ &&2\cdot I_2 \\ &&&\ddots \\ &&&& N\cdot I_2\endpmatrix,
\end{equation}
so that (\ref{expu1}) continues formally to hold, even though now $u$ is no more the solution of (\ref{wave}).\footnote{Observe that, for sake of brevity, we continue to use the same notation $\bfw,\bfq$, and $D$ used for the infinite expansion, though, hereafter, they will refer to the finite counterparts (\ref{wqD1}).}
Nevertheless, in the spirit of Fourier-Galerkin methods \cite{Boyd2001}, by requiring the residual obtained by plugging $u$ into (\ref{wave}) be orthogonal to the functional subspace\footnote{The same Fourier-Galerkin procedure will be used for the Hamiltonian PDEs studied in the following sections.}
 $${\cal V}_N = \mathrm{span}\left\{ c_0(x), s_1(x), c_1(x), \dots, s_N(x), c_N(x)\right\},$$ which, for fixed $t$, contains $u$, one obtains a finite set of $2N+1$ ODEs, formally still given by (\ref{wave1}), for which Theorem~\ref{waveth2} continues formally to hold, with the only exception that now $H(\bfq,\bfp)$ is no more equivalent to the Hamiltonian functional (\ref{waveH}), but only yields an approximation to it. Nevertheless, it is well known that, under suitable regularity assumptions on $f$ and the initial data $u_0$ and $v_0$, this truncated version converges exponentially to the original functional (\ref{waveH}), as $N\rightarrow\infty$  (this phenomenon is usually referred to as {\em spectral accuracy}), as well as the truncated version of $u$ converges to the infinite expansion (\ref{expu}). 

The resulting finite-dimensional semi-discrete problem (\ref{wave1}), which is still Hamiltonian, is the one we will solve by using line-integral methods. Actually, it is not yet ready to be solved, since the integral $\int_a^b \bfw(x)f'(\bfw(x)^\top\bfq)\dd x$, appearing in it, needs to be evaluated. For this purpose, since we are dealing with an integrand which is periodic in space, a composite trapezoidal rule based at the abscissae%\marginpar{\blue xim}
\begin{equation}\label{xim}
x_i = a + i\frac{b-a}m, \qquad i=0,\dots,m,
\end{equation}
can be considered. We refer, e.g., to \cite[Theorems\,5,\,6]{BFCI2015}, for a proper choice of the number, $m+1$, of points in (\ref{xim}), able to preserve the property of spectral accuracy.

Problem (\ref{wave1}) can then be solved by using a HBVM$(k,s)$ method, for which the accuracy results of Theorem~\ref{hbvmth} hold true. In particular, concerning the conservation of the semi-discrete Hamiltonian (\ref{waveH1}), the next result holds true, which follows from (\ref{Hcons}).%\marginpar{\red hbvmwave}

\begin{Theorem}\label{hbvmwave}
If a HBVM$(k,s)$ method is used with time-step $h$ for solving (\ref{wave1}), one has%\marginpar{\blue hbvmwave1}
\begin{equation}\label{hbvmwave1}
H(\bfq_1,\bfp_1)-H(\bfq_0,\bfp_0) = \left\{ \begin{array}{cc} 0, &\mbox{if}\quad f\in\Pi_\nu \quad\mbox{with}\quad \nu\le 2k/s,\\[2mm]
O(h^{2k+1}), &\mbox{otherwise}.\end{array}\right.
\end{equation}
having set $\bfq_1\approx \bfq(h)$ and $\bfp_1\approx \bfp(h)$ the new approximations.
\end{Theorem} 

\subsection{The nonlinear iteration}
In light of what previously stated, in order to obtain a {\em spectral accuracy in space} a suitably large value of $N$ in (\ref{wqD1}) has to be considered (a practical criterion for its choice will be sketched in Section~\ref{numsec}). Consequently, the special second-order problem (\ref{wave1}) is semilinear, with a bounded nonlinear term,\footnote{When the solution is bounded.} and the linear term given by $-D^2\bfq$. On the other hand, both the size (i.e., $2N+1$) and the norm \big(i.e., $\left(\frac{2\pi N}{b-a}\right)^2$\big) of the matrix $D^2$  tend to infinity, as $N\rightarrow\infty$. Consequently, when using a HBVM$(k,s)$ method for solving (\ref{wave1}), the blended iteration (\ref{ros}), (\ref{blend5})-(\ref{tetas1}) can be conveniently used, to get rid of the large norm of the linear term, with matrix $\Sigma$ approximated as in (\ref{Sigma2}). As a result, it turns out to be given by 
\begin{equation}\label{Sigmawave}
\Sigma = I_{2N+1} +h^2\rho_s^2 D^2,
\end{equation}
which is a {\em diagonal matrix} and, therefore, $\Sigma^{-1}$ can be cheaply computed and stored. Consequently, the complexity of the blended iteration turns out to be comparable with that of an explicit method, though not suffering from the step-size restrictions of this latter. As a matter of fact, the use of an explicit method usually would require $h\|D\|<1$, i.e., $h=O(N^{-1})$, which may be restrictive, when $N\gg1$. 

\subsection{Extension to higher space dimensions}\label{hisp}
For completeness, in this section we sketch the generalization of most of the previous arguments to the case where the space domain of the wave equation is, for sake of simplicity, the square $[a,b]^2:=[a,b]\times[a,b]$:
\begin{eqnarray}\label{wave2D}
u_{tt}(x,y,t) &=& \Delta u(x,y,t) - f'(u(x,t)), \qquad (x,y,t)\in[a,b]^2\times[0,T], \\
u(x,y,0)&=&u_0(x,y), \quad u_t(x,y,0)~=~v_0(x,y), \qquad (x,y)\in[a,b]^2.\nonumber\end{eqnarray}
As before, the problem (\ref{wave2D}) is completed by prescribing periodic boundary conditions. In such a case, the Hamiltonian functional becomes, by setting as usual\, $v=u_t$,
\begin{equation}\label{waveH2D}
\H[u,v](t) = \frac{1}2\int_a^b\int_a^b \left[ v^2(x,y,t)+\|\nabla u(x,y,t)\|_2^2 +2f(u(x,y,t))\right]\dd x\dd y.
\end{equation} 
Because of the periodic boundary conditions, we can again consider the orthonormal basis (\ref{c0})--(\ref{sj}) in each space dimension, thus obtaining the expansion
(for sake of brevity, let us set $s_0\equiv0$)
\begin{equation}\label{expu2D}
u(x,y,t) = \sum_{j,k\ge 0} \left[c_j(x)\alpha_j(t)+s_j(x)\beta_j(t)\right]\cdot\left[c_k(y)\eta_k(t)+s_k(y)\mu_k(t)\right],
\end{equation}
involving the additional time-dependent coefficients $\eta_0(t),\mu_1(t),\eta_1(t),\dots$. With reference to the infinite-dimensional vectors in (\ref{wq}), and defining the vectors
$$
\bfq_1(t):=\bfq(t), \qquad \bfq_2(t) = \pmatrix{cccc} \eta_0(t), &\mu_1(t), &\eta_1(t), &\dots\endpmatrix^\top,
$$
the infinite expansion (\ref{expu2D}) can be cast in vector form as
\begin{equation}\label{expu12D}
u(x,y,t) = [\bfw(x)\otimes\bfw(y)]^\top\,\bfq_1(t)\otimes\bfq_2(t).
\end{equation}
Consequently, by taking into account (\ref{Dmat})-(\ref{ortoabmat}), one obtains that (compare with Theorem~\ref{waveth2}) problem (\ref{wave2D}) can be recast as the infinite set of second order ODEs:
\begin{eqnarray}\nonumber
&&\!\!\!\!\!\!\!\!\!\ddot \bfq_1\otimes\ddot \bfq_2 =\! -D^2\bfq_1\otimes D^2 \bfq_2\! -\!\int_a^b\!\!\int_a^b\!\! \bfw(x)\otimes \bfw(y) f'\left([\bfw(x)\otimes\bfw(y)]^\top\,\!\!\bfq_1\otimes\bfq_2\right)\!\dd x\dd y, \, t\in[0,T], \\ 
\label{wave12D}%\\ \nonumber
&&\!\!\!\!\!\!\!\!\!\bfq_1(0)\otimes\bfq_2(0) =\int_a^b \int_a^b\bfw(x)\otimes\bfw(y)u_0(x,y)\dd x\dd y,\\ \nonumber
&&\!\!\!\!\!\!\!\!\!\dot\bfq_1(0)\otimes\dot\bfq_2(0) =\int_a^b\int_a^b \bfw(x)\otimes\bfw(y)v_0(x,y)\dd x\dd y.
\end{eqnarray}
By setting $\bfp_1\otimes\bfp_2=\dot\bfq_1\otimes\dot\bfq_2$, one then obtains that problem (\ref{wave12D}) is Hamiltonian, with Hamiltonian function
\begin{eqnarray}\nonumber
H(\bfq_1\otimes\bfq_2,\bfp_1\otimes\bfp_2) &=& \frac{1}2\left( (\bfp_1\otimes\bfp_2)^\top (\bfp_1\otimes\bfp_2) +
(\bfq_1\otimes\bfq_2)^\top (D\otimes D)^2 (\bfq_1\otimes\bfq_2)  \right.\\ \label{waveH12D}
&&\left.+ 2\int_a^b\int_a^b f\left([\bfw(x)\otimes\bfw(y)]^\top\,\bfq_1\otimes\bfq_2\right)\dd x\dd y\right).
\end{eqnarray}
This latter function, in turn, is equivalent to the Hamiltonian functional (\ref{waveH2D}), via the expansion (\ref{expu12D}). Then, as done in the one dimensional case, the vectors $\bfw(x), \bfw(y), \bfq_1(t), \bfq_2(t)$, are truncated after $2N+1$ terms, for a convenient large value of $N$, so that (\ref{wave12D}) becomes a Hamiltonian set of $(2N+1)^2$ ODEs, with Hamiltonian (\ref{waveH12D}).  This problem can be solved by adapting the arguments previously explained in the one dimensional case, even though now the complexity is clearly increased. Remarkably enough, however, the diagonal structure of the Jacobian of the linear term in (\ref{wave12D}), i.e., $-(D\otimes D)^2$, is still preserved.\footnote{Evidently, this property holds true whichever is the dimension of the considered space domain.}

%%%%%%%%%%%%%%%%%%%%%%%%%%%%%%%%%%%%%%%%%%
\section{The nonlinear Schr\"odinger equation}\label{NLSEsec}

We now consider the nonlinear Schr\"odinger equation, which is very important in many applications (see, e.g., the introduction in \cite{BBFCI2018}). In real variables, it takes the form,%\marginpar{\blue NLSE}
\begin{eqnarray}\label{NLSE}
u_t&=&-v_{xx} -f'(u^2+v^2)v, \qquad u(x,0) = u_0(x),\\
v_t&=&u_{xx}+f'(u^2+v^2)u, \qquad v(x,0) = v_0(x),\qquad (x,t)\in[a,b]\times[0,T],
\nonumber
 \end{eqnarray}
$f'$ being the derivative of a suitably regular function $f$. The problem is completed with periodic boundary conditions and, hereafter, we shall assume the initial functions to be suitably regular (as periodic functions), in order to guarantee a suitably smooth solution. Such an equation can be written in the form (\ref{J2}), with $\nabla\H$ the vector of the functional derivatives of the Hamiltonian functional%\marginpar{\blue NLSEH}
 \begin{equation}\label{NLSEH}
 \H[u,v](t) = \frac{1}2\int_a^b\left[ u_x^2 + v_x^2 -f(u^2+v^2) \right]\dd x.
 \end{equation}
This latter functional is conserved, because of the periodic boundary conditions \cite[Theorem\,1]{BBFCI2018}. Additional conserved (quadratic) functionals are the {\em mass} and the {\em momentum} \cite[Theorem\,2]{BBFCI2018}, respectively given by:%\marginpar{\blue NLSEMM} 
\begin{equation}\label{NLSEMM}
\M_1[u,v](t) = \int_a^b ( u^2+v^2 )\dd x, \qquad \M_2[u,v](t) = \frac{1}2\int_a^b ( v_xu-u_xv )\dd x.
\end{equation}
 
In order to obtain a space discretization which takes into account of the periodic boundary conditions, we consider again the expansion along the Fourier basis (\ref{c0})--(\ref{ortoab}), for $u$ and $v$. The expansion for $u$ is formally still given by (\ref{expu}). Similarly, that for $v$ will be given by:%\marginpar{\blue expv}  
\begin{equation}\label{expv}
v(x,t) = c_0(x)\eta_0(t)+\sum_{j\ge 1} \left[c_j(x)\eta_j(t)+s_j(x)\mu_j(t)\right],
\end{equation}
for suitable time dependent coefficients, $\eta_0(t), \eta_1(t),\mu_1(t),\dots$. By using the infinite vectors (\ref{wq}) and%\marginpar{\blue p}
\begin{equation}\label{p}
\bfp(t) = \pmatrix{cccc} \eta_0(t), &\mu_1(t), &\eta_1(t), &\dots\endpmatrix^\top,
\end{equation}
we can cast the expansions of $u$ and $v$ in vector form, respectively, as (\ref{expu1}) and%\marginpar{\blue expv1}
\begin{equation}\label{expv1}
v(x,t) = \bfw(x)^\top\bfp(t).
\end{equation}
As a consequence, the following result holds true, whose proof is similar to that of Theorem~\ref{waveth} (see also \cite[Section\,2]{BBFCI2018}).%\marginpar{\blue NLSEth}

\begin{Theorem}\label{NLSEth} With reference to (\ref{c0})--(\ref{ortoabmat}), problem (\ref{NLSE}) can be rewritten as the infinite-dimensional Hamiltonian ODE problem%\marginpar{\blue NLSE1}
\begin{eqnarray}\nonumber
\!\!\!\!\dot \bfq &=& D^2\bfp -\int_a^b \left[\bfw(x)f'((\bfw(x)^\top\bfq)^2+(\bfw(x)^\top\bfp)^2)\bfw(x)^\top\bfp \right]\dd x,  \\ \label{NLSE1}
\!\!\!\!\dot \bfp &=& -D^2\bfq +\int_a^b \left[\bfw(x)f'((\bfw(x)^\top\bfq)^2+(\bfw(x)^\top\bfp)^2)\bfw(x)^\top\bfq \right]\dd x, \,\, t\in[0,T], \\ \nonumber
\!\!\!\!\bfq(0) &=&\int_a^b \bfw(x)u_0(x)\dd x~=:~\bfq_0, \qquad \bfp(0) ~=~\int_a^b \bfw(x)v_0(x)\dd x~=:~\bfp_0.
\end{eqnarray}
This latter problem is Hamiltonian w.r.t. the Hamiltonian%\marginpar{\blue NLSEH1}
\begin{equation}\label{NLSEH1}
H(\bfq,\bfp) = \frac{1}2\left( \bfp^\top D^2 \bfp + \bfq^\top D^2\bfq - \int_a^b f((\bfw(x)^\top\bfq)^2+(\bfw(x)^\top\bfp)^2)\dd x\right),
\end{equation}
which turns out to be equivalent to the Hamiltonian functional (\ref{NLSEH}). Moreover, the two quadratic invariants (\ref{NLSEMM}) can be respectively rewritten as%\marginpar{\blue NLSEMM1}
\begin{equation}\label{NLSEMM1}
M_1(\bfq,\bfp) = \bfq^\top\bfq + \bfp^\top\bfp, \qquad M_2(\bfq,\bfp) = \bfq^\top \bar D\bfp,
\end{equation}
where $\bar D$ is the matrix defined in (\ref{bDmat}).
\end{Theorem}

As in the case of the nonlinear wave equation, in order to solve problem (\ref{NLSE1}) on a computer, one needs to truncate the infinite expansions (\ref{expu}) and (\ref{expv}) at a convenient index $N$. In so doing, the infinite vectors and matrices (\ref{wq}), (\ref{Dmat}), and (\ref{p}) become those in (\ref{wqD1}) and%\marginpar{\blue p1}
\begin{equation}\label{p1}
\bfp(t) = \pmatrix{cccccc} \eta_0(t), &\mu_1(t), &\eta_1(t), &\dots, &\mu_N(t), &\eta_N(t)\endpmatrix^\top,
\end{equation}
respectively. As a result, one eventually arrives again at the finite-dimensional Hamiltonian ODE problem (\ref{NLSE1}), having dimension $4N+2$, with the Hamiltonian and the invariants still given by (\ref{NLSEH1}) and (\ref{NLSEMM1}), respectively. Again, spectral accuracy is expected, if the solution is regular enough in space (as a periodic function). Finally, we mention that also in this case  the integrals in space can be computed by means of a composite trapezoidal rule, based at the abscissae (\ref{xim}), for a suitably large value of $m$.

Again, we can use an HBVM$(k,s)$ method for solving (\ref{NLSE1}). Concerning the conservation of the Hamiltonian, the following straightforward result follows from (\ref{Hcons}).%\marginpar{\red hbvmNLSE}

\begin{Theorem}\label{hbvmNLSE}
If a HBVM$(k,s)$ method is used with time-step $h$ for solving (\ref{NLSE1}), one has%\marginpar{\blue hbvmNLSE1}
\begin{equation}\label{hbvmNLSE1}
H(\bfq_1,\bfp_1)-H(\bfq_0,\bfp_0) = \left\{ \begin{array}{cc} 0, &\mbox{if}\quad f\in\Pi_\nu \quad\mbox{with}\quad \nu\le k/s,\\[2mm]
O(h^{2k+1}), &\mbox{otherwise}.\end{array}\right.
\end{equation}
having set $\bfq_1\approx \bfq(h)$ and $\bfp_1\approx \bfp(h)$ the new approximations.
\end{Theorem} 

\subsection{The nonlinear iteration}

Following the same arguments discussed in the previous section, in order to obtain a {\em spectral accuracy in space} a suitably large value of $N$ in (\ref{wqD1}) and (\ref{p1}) has to be considered (we remind that a practical criterion for its choice will be given in Section~\ref{numsec}). Consequently, the Hamiltonian problem (\ref{NLSE1}) is semilinear, with a bounded nonlinear term, if the solution is bounded, and the linear term given by (see (\ref{J2}))%%\marginpar{\blue NLSEJapp}
$$%\begin{equation}\label{\blue NLSEJapp}
J_2\otimes D^2\pmatrix{c}\bfq\\ \bfp\endpmatrix.
$$%\end{equation} 
On the other hand, the norm of the matrix $D^2$  is $\left(\frac{2\pi N}{b-a}\right)^2$, and tends to infinity, as $N\rightarrow\infty$, as well as its size. Consequently, when using a HBVM$(k,s)$ method for solving (\ref{NLSE1}), the blended iteration (\ref{tetas})--(\ref{ros}) can be conveniently used, to get rid of the large norm of the linear term, with matrix $\Sigma$ approximated as in (\ref{Sigma1})  and, in the present context, given by%\marginpar{\blue rimuovere $h$}
\begin{equation}\label{SigmaNLSE}
\Sigma = \pmatrix{cc} I_{2N+1} & - B\\ B& I_{2N+1}\endpmatrix, \qquad  B=h\rho_s D^2,
\end{equation}
which is a block matrix with diagonal blocks and, therefore, $\Sigma^{-1}$ can be cheaply computed and stored. As matter of fact, one has \cite[Theorem\,5]{BBFCI2018}:%%\marginpar{\blue Sigma1NLSE}
$$%\begin{equation}\label{Sigma1NLSE}
\Sigma^{-1} = \pmatrix{cc} \Gamma & B\cdot\Gamma \\ -B\cdot\Gamma &\Gamma\endpmatrix, \qquad \Gamma = (I_{2N+1}+B^2)^{-1}, 
$$%\end{equation}
which is again a block matrix with diagonal blocks (actually, two vectors are enough to store it). As a consequence, also in the present case the complexity of the blended iteration turns out to be comparable with that of an explicit method, though not suffering from step-size restrictions. As a matter of fact, the use of an explicit method would require $h\|D^2\|<1$, i.e., $h=O(N^{-2})$, which may be very restrictive, when $N\gg1$. 

%%%%%%%%%%%%%%%%%%%%%%%%%%%%%%%%%%%%%%%%%%
\section{The Korteweg--de Vries (KdV) equation}\label{KdVsec}
The last Hamiltonian PDE that we consider is the Korteweg--de Vries equation, recently investigated in \cite{BGS2019} by using line integral methods,
\begin{eqnarray}\label{KdV}
u_t &=& \alpha u_{xxx} + \beta uu_x, \qquad (x,t)\in[a,b]\times[0,T], \\
u(x,0) &=& u_0(x),\nonumber 
\end{eqnarray}
with $\alpha\beta\ne0$, and coupled with periodic boundary conditions. As usual, we shall assume that $u_0$ is smooth enough, as a periodic function, so that $u(x,t)$ turns out to be suitably regular, as a periodic function in space \cite{KP2008}. Equation (\ref{KdV}) can be written in Hamiltonian form as
$$u_t = \partial_x \left(\delta_u \H[u]\right)$$
with $\H[u]$ the Hamiltonian functional%\marginpar{\blue KdVH}
\begin{equation}\label{KdVH}
\H[u](t) = \frac{1}2\int_a^b\left[ -\alpha u_x(x,t)^2+\frac{\beta}3 u(x,t)^3\right]\dd t \equiv \int_a^b L(x,t,u,u_x)\dd x,
\end{equation}
and
$$\delta_u\H[u] =  (\partial_u -\partial_x\partial_{u_x}) L(x,t,u,u_x),$$
its functional derivative.\footnote{Actually, it can be seen that there is a further Hamiltonian formulation of (\ref{KdV}) \cite{O1984}, so that the PDE has a so called {\em bi-Hamiltonian} structure.} Because of the periodic boundary conditions, the Hamiltonian functional (\ref{KdVH}) turns out to be conserved. Another conserved functional is given by%\marginpar{\blue KdVU}
\begin{equation}\label{KdVU}
\U[u] = \int_a^b u\,\dd x,
\end{equation}
as it can be readily shown. In order to obtain a space discretization, we consider an expansion along the usual orthonormal basis (\ref{c0})--(\ref{ortoab}), which provides us with an expression formally still given by (\ref{expu}). However, because of the conservation of the functional (\ref{KdVU}), one obtains that
$$
\int_a^b u_0(x)\dd x = \int_a^b u(x,t)\dd x \equiv \int_a^b c_0(x)\alpha_0(t)\dd x = (b-a)c_0(x)\alpha_0(t), \quad \forall t\ge0.
$$
Consequently, the expansion (\ref{expu}) now becomes%\marginpar{\blue KdVexpu}
\begin{equation}\label{KdVexpu}
u(x,t) = \hat u_0+\sum_{j\ge 1} \left[c_j(x)\alpha_j(t)+s_j(x)\beta_j(t)\right], \qquad \hat u_0 = (b-a)^{-1}\int_a^bu_0(x)\dd x.
\end{equation}
In order to put this expansion in vector form, let us introduce the infinite vectors%\marginpar{\blue csqp}
\begin{equation}\label{csqp}
\bfc(x) = \pmatrix{c} c_1(x)\\ c_2(x)\\ \vdots\endpmatrix,\quad
\bfs(x) = \pmatrix{c} s_1(x)\\ s_2(x)\\ \vdots\endpmatrix,\quad
\bfq(t) = \pmatrix{c} \alpha_1(t)\\ \alpha_2(t)\\ \vdots\endpmatrix,\quad
\bfp(t) = \pmatrix{c} \beta_1(t)\\ \beta_2(t)\\ \vdots\endpmatrix.
\end{equation}
In so doing, we can rewrite (\ref{KdVexpu}) as:%\marginpar{\blue KdVexpu1}
\begin{equation}\label{KdVexpu1}
u(x,t) = \hat u_0+\bfc(x)^\top\bfq(t)+\bfs(x)^\top\bfp(t), \qquad \hat u_0 = (b-a)^{-1}\int_a^bu_0(x)\dd x.
\end{equation}
Consequently, the conservation of (\ref{KdVU}) is automatically granted. 
Moreover, by defining the infinite matrix%\marginpar{\blue KdVD}
\begin{equation}\label{KdVD}
D = \frac{2\pi}{b-a}\pmatrix{ccc} 1\\ &2\\ &&\ddots\endpmatrix,
\end{equation}
such that%\marginpar{\blue c1s1}  
\begin{equation}\label{c1s1}
\bfc'(x) = -D\bfs(x), \qquad \bfs'(x)=D\bfc(x),
\end{equation}  
and similarly for the higher derivatives, and considering that%\marginpar{\blue ccsscs}
\begin{equation}\label{ccsscs}
\int_a^b \bfc(x)\bfc(x)^\top\dd x = \int_a^b \bfs(x)\bfs(x)^\top\dd x = I, \qquad \int_a^b \bfc(x)\bfs(x)^\top\dd x=O,
\end{equation}
one verifies that (\ref{KdV}) can be rewritten as the infinite dimensional ODE problem (see \cite[Lemma\,3]{BGS2019} for full details)%\marginpar{\blue KdV1}
\begin{eqnarray}\nonumber
\dot\bfq &=& D\left[ -\alpha D^2\bfp + \frac{\beta}2\int_a^b \bfs\, (\hat u_0+\bfc^\top\bfq+\bfs^\top\bfp)^2\dd x\right],\\ \label{KdV1}
\dot\bfp &=& -D\left[ -\alpha D^2\bfq + \frac{\beta}2\int_a^b \bfc\, (\hat u_0+\bfc^\top\bfq+\bfs^\top\bfp)^2\dd x\right],	\qquad t\in[0,T],\\
\bfq(0) &=& \int_a^b \bfc(x) u_0(x)\dd x ~=:~\bfq_0, \qquad \bfp(0) ~=~ \int_a^b \bfs(x) u_0(x)\dd x ~=:~\bfp_0.\nonumber
\end{eqnarray}
For this problem, the following result holds true \cite[Theorem\,1]{BGS2019}.%\marginpar{\red KdVth}

\begin{Theorem}\label{KdVth}
Problem (\ref{KdV1}) is in the form (\ref{Hpro}) with (see (\ref{J2}))
$$y = \pmatrix{c} \bfq\\ \bfp\endpmatrix, \qquad J = J_2\otimes D,$$ and the Hamiltonian given by%\marginpar{\blue KdVH1}
\begin{equation}\label{KdVH1}
H(\bfq,\bfp) = \frac{1}2\left[ -\alpha\left( \bfq^\top D^2 \bfq + \bfp^\top D^2 \bfp\right) +\frac{\beta}3\int_a^b (\hat u_0+\bfc^\top\bfq+\bfs^\top\bfp)^3\dd x\right].
\end{equation}
This latter is equivalent to the Hamiltonian functional (\ref{KdVH}), via (\ref{KdVexpu1}) and (\ref{c1s1})-(\ref{ccsscs}).
\end{Theorem}

As done before, in order for the problem (\ref{KdV1}) to be solvable on a computer, the infinite expansion in (\ref{KdVexpu}) must be truncated to a convenient index $N$. In so doing, one still formally retrieves the vector formulation (\ref{KdVexpu1}), where now the vectors%\marginpar{\blue csqp1}
\begin{equation}\label{csqp1}
\bfc(x) = \pmatrix{c} c_1(x)\\ \vdots \\ c_N(x)\endpmatrix,\quad
\bfs(x) = \pmatrix{c} s_1(x)\\ \vdots \\  s_N(x)\endpmatrix,\quad
\bfq(t) = \pmatrix{c} \alpha_1(t)\\ \vdots\\ \alpha_N(t)\endpmatrix,\quad
\bfp(t) = \pmatrix{c} \beta_1(t)\\ \vdots\\ \beta_N(t)\endpmatrix,
\end{equation}
are hereafter used in place of (\ref{csqp}). Similarly, by replacing matrix (\ref{KdVD}) with%\marginpar{\blue KdVD1}
\begin{equation}\label{KdVD1}
D = \frac{2\pi}{b-a}\pmatrix{ccc} 1\\ &\ddots\\ && N\endpmatrix,
\end{equation}
one obtains a set of $2N$ Hamiltonian equations, formally still given by (\ref{KdV1}), with the Hamiltonian $H$ also formally given by (\ref{KdVH1}). As for the Hamiltonian PDEs previously studied, spectral accuracy is expected, as $N\rightarrow\infty$, upon regularity assumptions on $u_0$.\footnote{Concerning the integrals appearing in (\ref{KdV1}) and (\ref{KdVH1}), they can be exactly computed via a composite trapezoidal rule based at the abscissae (\ref{xim}), by choosing $m>3N$ \cite{BGS2019}.}

Having got the finite dimensional Hamiltonian ODE problem (\ref{KdV1}), we can use a HBVM$(k,s)$ method for its time integration. Concerning energy conservation, the following result easily follows from (\ref{Hcons}).%\marginpar{\red KdVth1}

\begin{Theorem}\label{KdVth1}
A HBVM$(k,s)$ method used for solving (\ref{KdV1}) is energy-conserving, for all $k\ge 3s/2$.
\end{Theorem}

\subsection{The nonlinear iteration}
Also in this case, problem (\ref{KdV1}) is semilinear. However,  it is worth observing that the Hessian of the Hamiltonian $H$ in (\ref{KdVH1}) is given by
$$\nabla^2 H(\bfq,\bfp) = \pmatrix{cc} -\alpha D^2+\beta\int_a^b u(x,t) \bfc(x)\bfc(x)^\top\dd x & \beta \int_a^b u(x,t)\bfc(x)\bfs(x)^\top\dd x\\
\beta \int_a^b u(x,t)\bfs(x)\bfc(x)^\top\dd x & -\alpha D^2+\beta\int_a^b u(x,t) \bfs(x)\bfs(x)^\top\dd x\endpmatrix,$$
with $u(x,t)$ given by the expansion (\ref{KdVexpu1}).
Consequently, by considering the constant approximation $u(x,t)\equiv \hat u_0$ (due to the conservation of (\ref{KdVU})), and taking into account (\ref{ccsscs}), one obtains the {\em constant approximate (diagonal) Hessian}
$$\nabla^2 H(\bfq,\bfp)\approx I_2\otimes \hat D, \qquad \hat D:=\left[ -\alpha D^2+\beta \hat u_0I_N\right].$$ Therefore, the blended iteration (\ref{tetas})--(\ref{ros}) can be conveniently used, by considering the resulting approximated matrix (see (\ref{J2}) and (\ref{KdVD1}))%\marginpar{\blue rimuovere $h$}
$$\Sigma = I_2\otimes I_N - h\rho_s(J_2\otimes D)(I_2\otimes \hat D) \equiv  
\pmatrix{cc} I_N & -B\\ B& I_N\endpmatrix, \qquad B=h\rho_s D \hat D,$$
which is a block matrix with diagonal blocks. Moreover, one has \cite[Theorem\,3]{BGS2019}:
$$
\Sigma^{-1} = \pmatrix{cc} \Gamma & B\cdot\Gamma \\ -B\cdot\Gamma &\Gamma\endpmatrix, \qquad \Gamma = (I_N+B^2)^{-1}, 
$$
which can be easily computed (once for all) and stored (in fact, only two vectors of length $N$ are needed). Consequently, the complexity of the blended iteration turns out to be comparable with that of an explicit method, though not suffering from its step-size restrictions which, for the present problem, would require  $h=O(N^{-3})$. 

%%%%%%%%%%%%%%%%%%%%%%%%%%%%%%%%%%%%%%%%%%
\section{Numerical tests}\label{numsec}

In this section, we report a few numerical tests, aimed at assessing the effectiveness of HBVMs for solving the previously studied Hamiltonian PDEs. In particular, the spectral version of HBVMs (SHBVMs) will be recognized to be very promising. In more details, we shall compare the following methods:
\begin{itemize}
\item the symplectic $s$-stage Gauss methods, $s=1,2$;

\item the energy-conserving HBVM$(k,s)$ methods, $s=1,2$, and $k$ suitably chosen;

\item the SHBVM method.
\end{itemize}
The comparisons will be quite fair, since the same Matlab function\,\footnote{It is a modification of the function {\tt hbvm} available at \cite{LIMbook}.}  implements all methods. All numerical tests have been done on a 2.8 GHz Intel Core i7 computer with 16GB of memory, running Matlab 2017b. 

To begin with, let us define the criterion used for getting spectral accuracy in space, i.e., for a correct choice of $N$ in (\ref{wqD1}), (\ref{p1}), (\ref{csqp1}), and (\ref{KdVD1}). In more details, $N$ has been chosen in order to fulfil both the two following requirements:
\begin{itemize}
\item {\em a good approximation of the initial condition.} This is achieved by requiring%\marginpar{\blue ini1}
\begin{equation}\label{ini1}
E_0 ~:=~\max\left\{\|u_0(x)-\bfw(x)^\top\bfq_0\|_\infty,\, \|v_0(x)-\bfw(x)^\top\bfp_0\|_\infty\right\} ~\le~ tol ~\approx~ \varepsilon,
\end{equation}
with $\varepsilon$ the machine epsilon, for problems (\ref{wave1}) and (\ref{NLSE1}), or%\marginpar{\blue ini2}
\begin{equation}\label{ini2}
E_0~:=~\|u_0(x)-\hat u_0-\bfc(x)^\top\bfq_0-\bfs(x)^\top\bfp_0\|_\infty ~\le~ tol~ \approx~ \varepsilon,
\end{equation}
for problem (\ref{KdV1});

\smallskip
\item {\em a good approximation of the Hamiltonian.} This is achieved by computing the initial value $H(\bfq_0,\bfp_0)=:H_0$ of the semi-discrete Hamiltonian (i.e., (\ref{waveH1}), or (\ref{NLSEH1}), or (\ref{KdVH1})) for consecutive values of $N$, and checking that the absolute value of the difference, $\Delta H_0$, satisfies:%\marginpar{\blue iniH}
\begin{equation}\label{iniH}
\Delta H_0 ~\le~tol~\approx~\varepsilon.
\end{equation} 

\end{itemize}

\subsection{The semilinear wave equation} 
\begin{figure}[t]
\centering
\includegraphics[width=9 cm,height=6 cm]{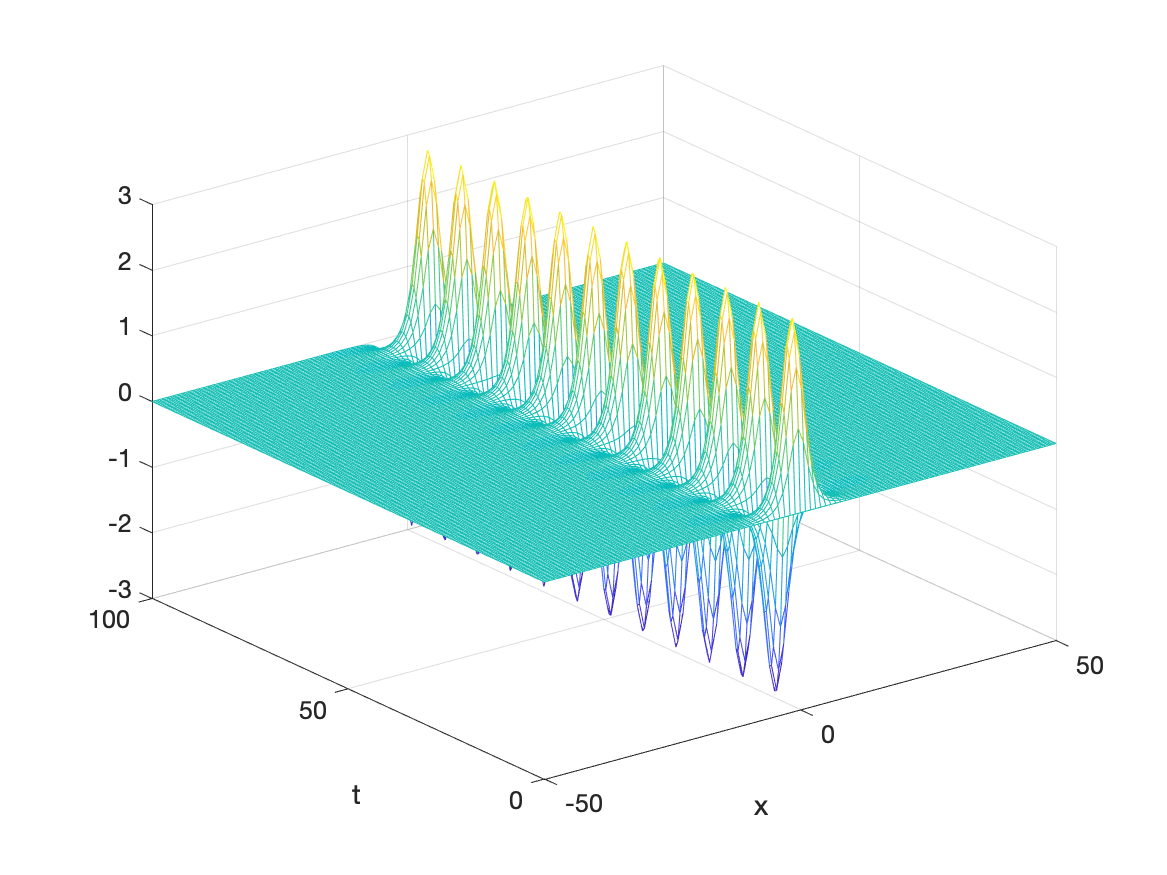} \qquad \includegraphics[width=9 cm,height=6 cm]{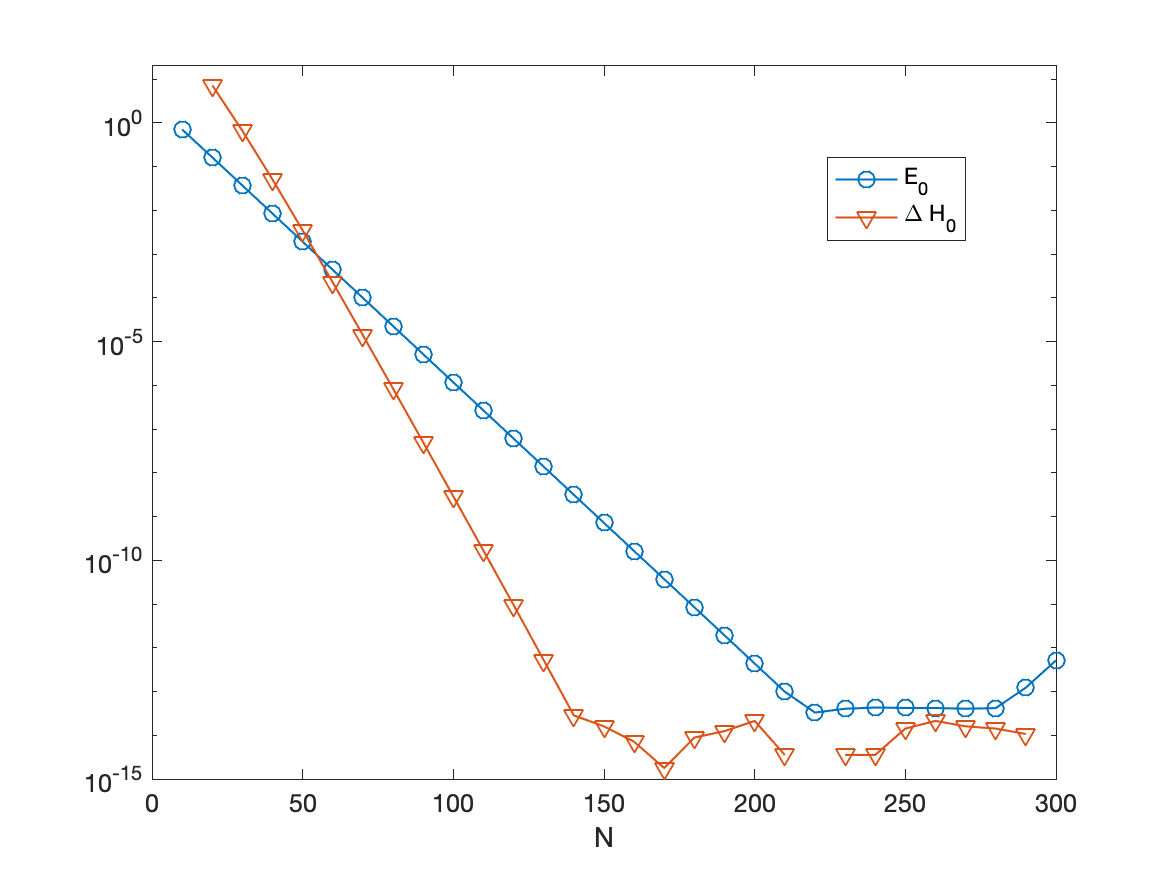}
\caption{Sine-Gordon equation (\ref{sGe}); upper plot: solution (\ref{sGeu}); lower plot: $E_0$ (see (\ref{ini1})) and $\Delta H_0$ (see (\ref{iniH})) versus $N$.}\label{sGfig}
\end{figure}

We consider the so called {\em sine-Gordon} equation \cite[Section\,7]{BFCI2015} with a {\em breather soliton} solution,%\marginpar{\blue sGe}
\begin{eqnarray}\label{sGe}
u_{tt}&=& u_{xx} -\sin(u), \qquad (x,t)\in[-50,50]\times[0,100],\\ \nonumber
u(x,0) &=&0, \qquad u_t(x,0) = \frac{4}\gamma\mathrm{sech}\left(\frac{x}\gamma\right),
\end{eqnarray}
where we choose $\gamma=1.5$. Its solution, depicted in the upper plot in Figure~\ref{sGfig}, is:%\marginpar{\blue sGeu}
\begin{equation}\label{sGeu}
u(x,t) = 4\,\mathrm{atan}\hspace{-.3em}\left( \mathrm{sech}\left(\frac{x}\gamma\right) \frac{\sin\left( t\sqrt{1-\gamma^{-2}}\right)}{\sqrt{\gamma^2-1}}\right).
\end{equation}
In the lower plot of Figure~\ref{sGfig} there are the graphs of $E_0$ and $\Delta H_0$, as defined in (\ref{ini1}) and (\ref{iniH}), respectively. From such plots, one infers that the choice $N=250$ is adequate to obtain spectral accuracy in space. In Table~\ref{sGtab} we list the obtained numerical results by solving the resulting semi-discrete problem (\ref{wave1}) with time-step $h=100/n$. In more details: the execution time (in sec), the maximum solution and Hamiltonian errors, $e_u$ and $e_H$, respectively, and the rate of convergence, where appropriate; for the SHBVM method, we also list the used values of $k$ and $s$, the latter obtained by using $tol\approx \sqrt{\varepsilon}$ in (\ref{tol}) and $k$ suitably larger than $s$. From the obtained results, one sees that:
\begin{itemize}
\item the higher-order methods perform better than the lower-order ones;
\item the energy-conserving methods are slightly more efficient than the symplectic ones, when the largest time-steps are used; 
\item the spectral method turns out to be the most effective one, and uses much larger time-steps.
\end{itemize}

\begin{table}[t]
\caption{Numerical solution of the sine-Gordon equation (\ref{sGe}) using a time-step $h=100/n$.}\label{sGtab}
\centering
\begin{tabular}{rccccc}
\toprule
\multicolumn{6}{c}{Gauss 1}\\
\toprule
$n$~ & time & $e_u$ & rate & $e_H$ & rate\\
\midrule
2000 &  2.1 & 4.61e-02 & --- & 1.98e-03 & ---  \\ 
3000 &  2.8 & 2.05e-02 & 2.0 & 8.79e-04 & 2.0  \\ 
4000 &  3.8 & 1.15e-02 & 2.0 & 4.95e-04 & 2.0  \\ 
5000 &  4.8 & 7.37e-03 & 2.0 & 3.17e-04 & 2.0  \\ 
6000 &  6.1 & 5.12e-03 & 2.0 & 2.20e-04 & 2.0  \\
\bottomrule
\toprule
\multicolumn{6}{c}{Gauss 2}\\
\toprule
$n$~ & time & $e_u$ & rate & $e_H$ & rate\\
\midrule
1000 &  2.4 & 2.69e-05 & --- & 2.57e-06 & ---  \\ 
1500 &  3.3 & 5.28e-06 & 4.0 & 5.05e-07 & 4.0  \\ 
2000 &  3.9 & 1.67e-06 & 4.0 & 1.59e-07 & 4.0  \\ 
2500 &  5.3 & 6.83e-07 & 4.0 & 6.53e-08 & 4.0  \\ 
3000 &  6.6 & 3.29e-07 & 4.0 & 3.15e-08 & 4.0  \\ 
\bottomrule
\toprule
\multicolumn{6}{c}{HBVM(4,1)}\\
\toprule
$n$~ & time & $e_u$ & rate & $e_H$ \\
\midrule
1000 &  2.4 & 1.37e-02 & --- & 7.11e-15  \\ 
1500 &  3.3 & 6.15e-03 & 2.0 & 1.07e-14  \\ 
2000 &  4.2 & 3.47e-03 & 2.0 & 8.88e-15  \\ 
2500 &  5.7 & 2.22e-03 & 2.0 & 1.07e-14  \\ 
3000 &  7.0 & 1.55e-03 & 2.0 & 8.88e-15  \\ 
\bottomrule
\toprule
\multicolumn{6}{c}{HBVM(4,2)}\\
\toprule
$n$~ & time & $e_u$ & rate & $e_H$ \\
\midrule
1000 &  2.8 & 2.11e-05 & --- & 1.07e-14  \\ 
1500 &  4.0 & 4.18e-06 & 4.0 & 8.88e-15  \\ 
2000 &  4.6 & 1.32e-06 & 4.0 & 1.07e-14  \\ 
2500 &  6.0 & 5.42e-07 & 4.0 & 7.11e-15  \\ 
3000 &  7.0 & 2.61e-07 & 4.0 & 8.88e-15  \\ 
\bottomrule
\multicolumn{6}{c}{SHBVM}\\
\toprule
$n$~ & time & $k$ & $s$ & $e_u$ & $e_H$ \\
\midrule
 50 &  2.7 &  22 &  20 & 2.87e-12 & 3.55e-15  \\ 
 75 &  1.6 &  20 &  18 & 3.61e-13 & 7.11e-15  \\ 
100 &  1.3 &  15 &  12 & 3.53e-13 & 3.55e-15  \\ 
\bottomrule
\end{tabular}
\end{table}

\subsection{The nonlinear Schr\"odinger equation} 
\begin{figure}[t]
\centering
\includegraphics[width=9 cm,height=6 cm]{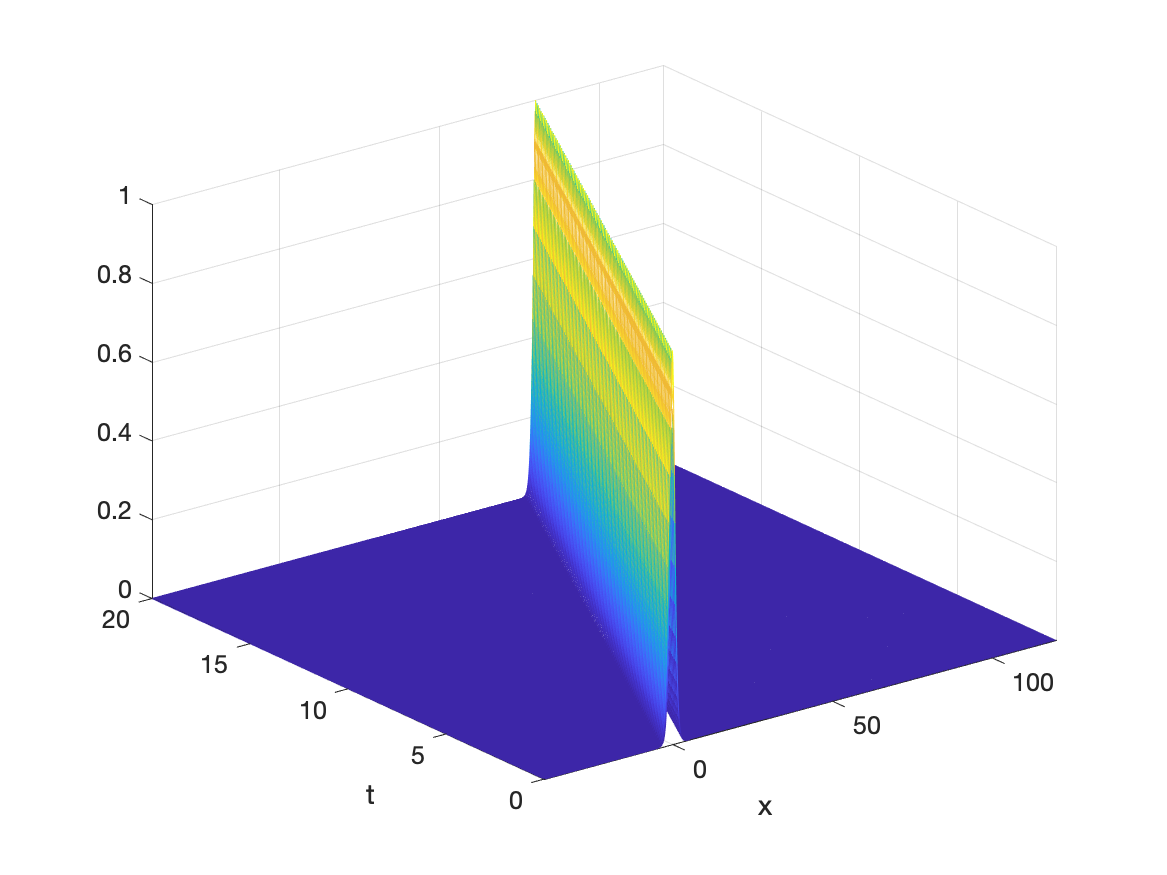} \qquad \includegraphics[width=9 cm,height=6 cm]{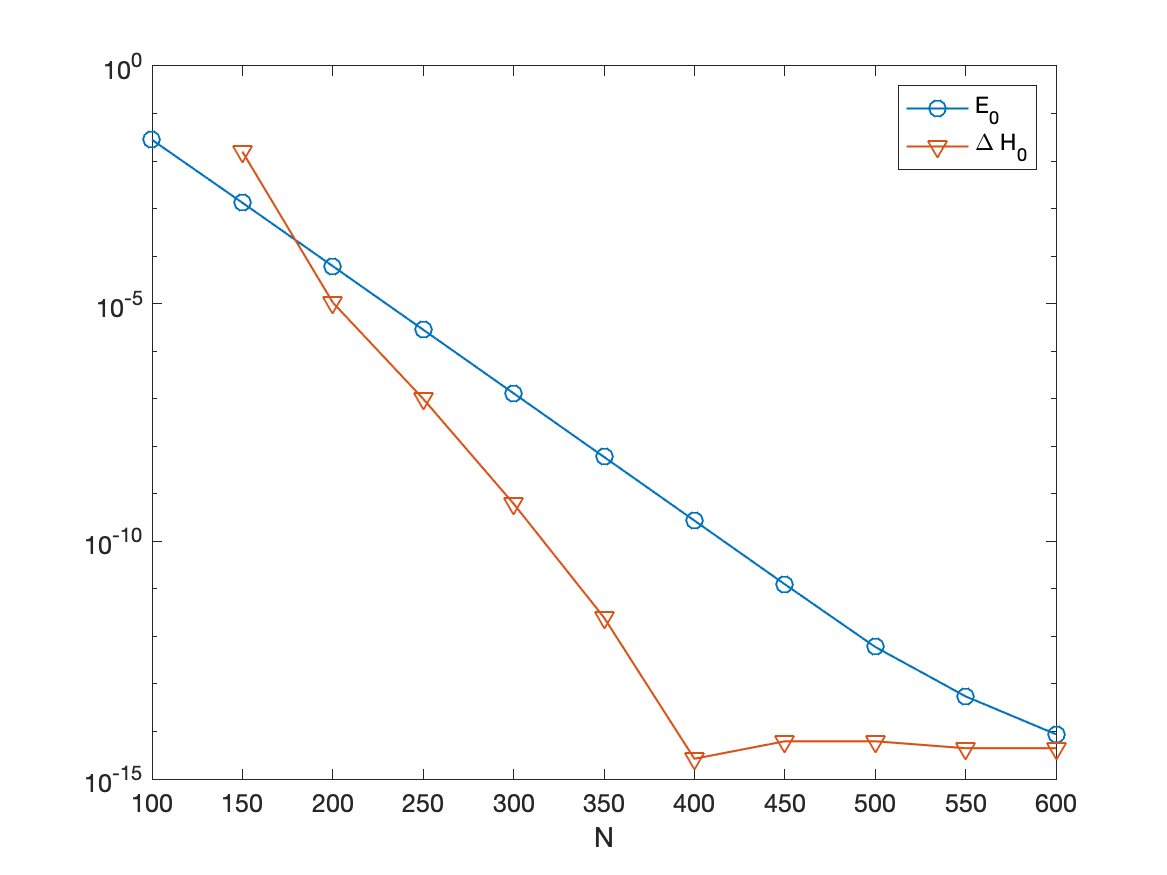}
\caption{Nonlinear Schr\"odinger equation (\ref{NLSEex}); upper plot: modulus of the solution (\ref{NLSEuv}); lower plot: $E_0$ (see (\ref{ini1})) and $\Delta H_0$ (see (\ref{iniH})) versus $N$.}\label{NLSEfig}
\end{figure}

We consider the so called {\em focusing} equation,\footnote{The {\em de-focusing} case is obtained when the sign of the coupling term is reversed.}%\marginpar{\blue NLSEex}
\begin{eqnarray}\nonumber
u_t&=& -v_{xx} -2(u^2+v^2)v, \\ \label{NLSEex} 
v_t&=& u_{xx} +2(u^2+v^2)u, \qquad (x,t)\in[-40,120]\times[0,20],
\end{eqnarray}
where the initial conditions at $t=0$ are taken from the known solution,%\marginpar{\blue NLSEuv}
\begin{equation}\label{NLSEuv}
u(x,t) =   \mathrm{sech}(x-4t)\cos(2x-3t), \qquad v(x,t) = \mathrm{sech}(x-4t)\sin(2x-3t),
\end{equation}
depicted in the upper plot of Figure~\ref{NLSEfig}, plus (approximate) boundary conditions. 
In the lower plot of the same figure, there are the plots of $E_0$ and $\Delta H_0$, as defined in (\ref{ini1}) and (\ref{iniH}), respectively. From such plots, one infers that the choice $N=600$ is adequate to obtain spectral accuracy in space.
For this problem, the symplectic $s$-stage Gauss methods conserve the quadratic invariants (\ref{NLSEMM1}), whereas the HBVM$(2s,s)$ methods are energy conserving (according to Theorem~\ref{hbvmNLSE}, since $f(x)=x^2$). For the SHBVM method, we use $tol\approx 10^{-1}\sqrt{\varepsilon}$ in (\ref{tol}). 
 In Table~\ref{NLSEtab} we list the numerical results  obtained by solving the resulting semi-discrete problem (\ref{NLSE1}) with time-step $h=20/n$: besides the execution time (in sec), we list the maximum solution, mass, momentum,  and Hamiltonian errors, $e_{uv}$, $e_1$, $e_2$, and $e_H$, respectively, along with the rate of convergence, where appropriate; for the SHBVM method, we also list the used values of $k$ and $s$. We observe that a kind of super-convergence occurs in the invariants (twice the convergence order of the solution) for the Gauss and HBVM methods. In this case, the symplectic and energy-conserving methods turn out to be almost equivalent, with the higher-order methods more efficient than the lower-order ones. However, the SHBVM method outperform all of them, being able to use much larger time-steps, and having a uniformly small error in both the solution and the invariants (which are all conserved within the round-off error level). 

\begin{table}[t]
\caption{Numerical solution of the nonlinear Schr\"odinger equation (\ref{NLSEex}) using a time-step $h=20/n$.}\label{NLSEtab}
\centering
\begin{tabular}{rcccccccc}
\toprule
\multicolumn{9}{c}{Gauss 1}\\
\toprule
$n$~ & time & $e_{uv}$ & rate & $e_1$ & $e_2$ & $e_H$ & rate\\
\midrule
400 & 19.1 & 4.93e-01 & --- & 1.60e-14 & 2.01e-16 & 9.58e-04 & ---  \\ 
600 & 26.3 & 2.32e-01 & 1.9 & 4.75e-14 & 6.70e-16 & 1.79e-04 & 4.1  \\ 
800 & 33.2 & 1.31e-01 & 2.0 & 2.38e-14 & 1.39e-16 & 5.55e-05 & 4.1  \\ 
1000 & 40.7 & 8.41e-02 & 2.0 & 2.22e-14 & 5.20e-17 & 2.25e-05 & 4.0  \\ 
\bottomrule
\toprule
\multicolumn{9}{c}{Gauss 2}\\
\toprule
$n$~ & time & $e_{uv}$ & rate & $e_1$ & $e_2$ & $e_H$ & rate\\
\midrule
400 & 41.3 & 1.71e-03 & --- & 3.13e-14 & 3.12e-17 & 5.20e-08 & ---  \\ 
600 & 60.4 & 3.41e-04 & 4.0 & 1.91e-14 & 2.43e-17 & 2.14e-09 & 7.9  \\ 
800 & 72.1 & 1.08e-04 & 4.0 & 1.91e-14 & 2.78e-17 & 2.18e-10 & 7.9  \\ 
1000 & 84.7 & 4.44e-05 & 4.0 & 2.00e-14 & 3.82e-17 & 3.69e-11 & 8.0  \\ 
\bottomrule
\toprule
\multicolumn{9}{c}{HBVM(2,1)}\\
\toprule
$n$~ & time & $e_{uv}$ & rate & $e_1$ & rate & $e_2$ & rate & $e_H$ \\
\midrule
400 & 36.3 & 5.23e-01 & --- & 1.40e-04 & --- & 4.25e-06 & --- & 3.55e-15  \\ 
600 & 50.7 & 2.45e-01 & 1.9 & 2.64e-05 & 4.1 & 7.96e-07 & 4.1 & 3.55e-15  \\ 
800 & 60.1 & 1.38e-01 & 2.0 & 8.23e-06 & 4.1 & 2.47e-07 & 4.1 & 4.00e-15  \\ 
1000 & 74.5 & 8.89e-02 & 2.0 & 3.35e-06 & 4.0 & 1.01e-07 & 4.0 & 4.00e-15  \\ 
\bottomrule
\toprule
\multicolumn{9}{c}{HBVM(4,2)}\\
\toprule
$n$~ & time & $e_{uv}$ & rate & $e_1$ & rate & $e_2$ & rate & $e_H$ \\
\midrule
400 & 43.2 & 1.74e-03 & --- & 6.66e-09 & --- & 1.83e-10 & --- & 4.44e-15  \\ 
600 & 61.4 & 3.47e-04 & 4.0 & 2.70e-10 & 7.9 & 7.48e-12 & 7.9 & 3.55e-15  \\ 
800 & 77.0 & 1.10e-04 & 4.0 & 2.74e-11 & 8.0 & 7.61e-13 & 7.9 & 4.00e-15  \\ 
1000 & 92.3 & 4.52e-05 & 4.0 & 4.63e-12 & 8.0 & 1.29e-13 & 8.0 & 3.55e-15  \\ 
\bottomrule
\multicolumn{9}{c}{SHBVM}\\
\toprule
$n$~ & time & $k$ & $s$ & $e_{uv}$ & $e_1$ & $e_2$ & $e_H$ \\
\midrule
 50 & 55.0 &  20 &  18 & 3.13e-11 & 1.35e-14 &6.59e-17 &4.88e-15  \\ 
 75 & 53.6 &  16 &  14 & 2.27e-11 & 1.33e-14 &7.29e-17 &3.11e-15  \\ 
100 & 61.6 &  14 &  12 & 2.47e-11 & 1.40e-14 &6.25e-17 &3.11e-15  \\ 
\bottomrule
\end{tabular}
\end{table}

\subsection{The Korteweg--de Vries equation}
\begin{figure}[t]
\centering
\includegraphics[width=9 cm,height=6 cm]{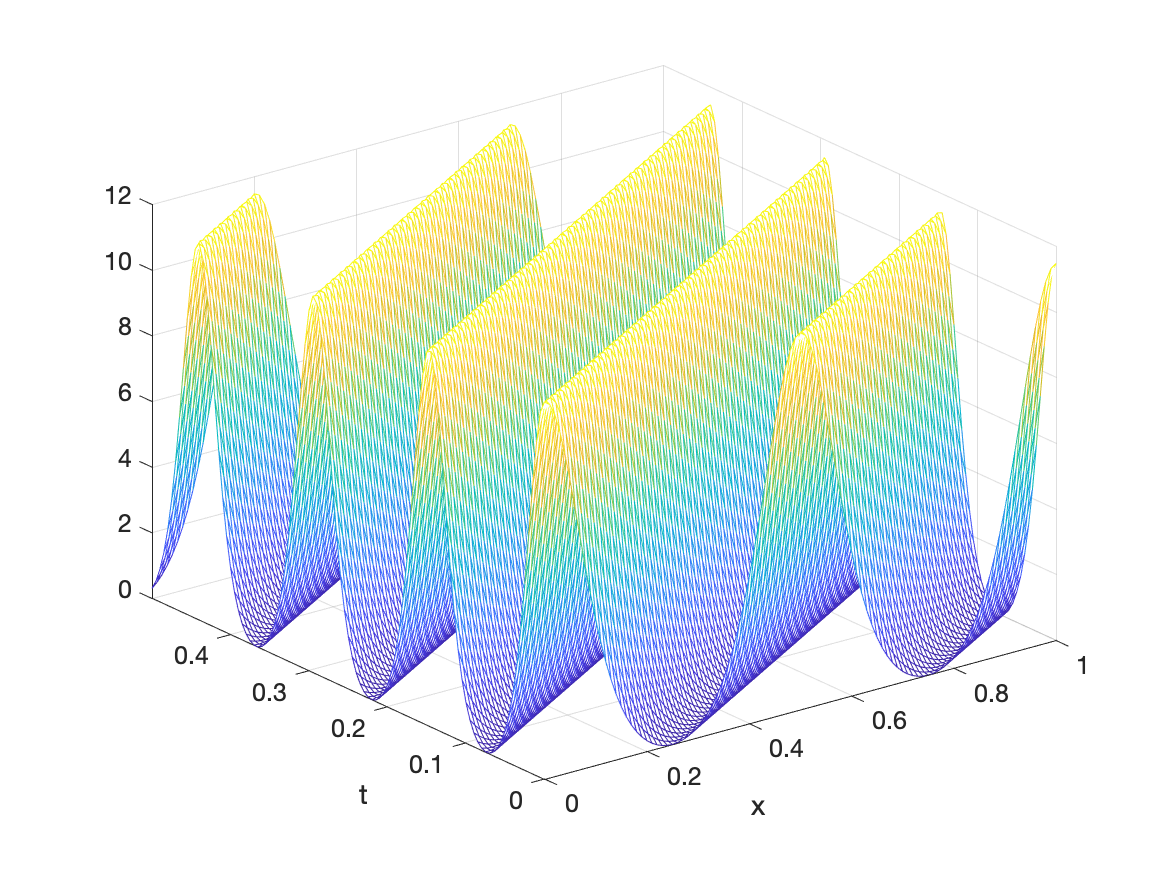} \qquad \includegraphics[width=9 cm,height=6 cm]{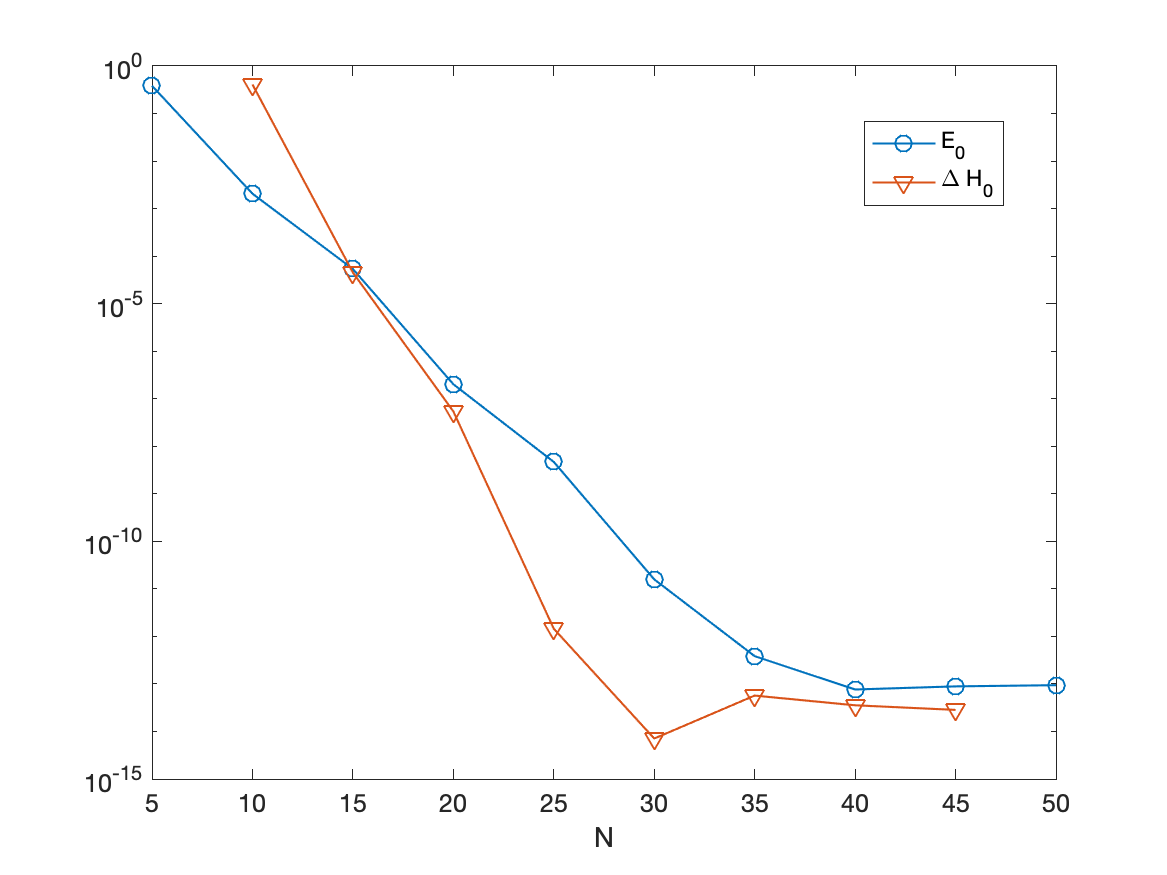}
\caption{Korteweg--de Vries equation (\ref{KdVe}); upper plot: solution (\ref{KdVeu}); lower plot: $E_0$ (see (\ref{ini2})) and $\Delta H_0$ (see (\ref{iniH})) versus $N$.}\label{KdVfig}
\end{figure}

This example is adapted from \cite[Example\,2]{BGS2019}:%\marginpar{\blue KdVe}
\begin{equation}\label{KdVe}
u_t+\epsilon u_{xxx} +u u_x = 0, \qquad (x,t)\in[0,1]\times[0,10],
\end{equation}
equipped with periodic boundary conditions and the initial condition obtained from the known {\em cnoidal wave solution},%\marginpar{\blue KdVeu}
\begin{equation}\label{KdVeu}
u(x,t) = a\, \mathrm{cn}^2\left(4K(m)(x-\nu t-x_0)\right).
\end{equation}
Here $\mathrm{cn}:=\mathrm{cn}(z|m)$ is the Jacobi elliptic function with modulus $m$, $K(m)$ is the complete elliptic integral of the first kind, and the following parameters have been used:
$$\epsilon = 10^{-2}, \qquad m = 0.9, \qquad a = 192 m\epsilon K^2(m), \qquad \nu=64\epsilon(2m-1)K^2(m), \qquad x_0=1/2.$$ 
The initial part of the solution (\ref{KdVeu}) is depicted in the upper plot of Figure~\ref{KdVfig}, whereas in the lower plot one may find $E_0$ and $\Delta H_0$, as defined in (\ref{ini2}) and (\ref{iniH}), respectively, versus $N$ . From the latter plots, one infers that the choice $N=50$ is adequate to obtain spectral accuracy in space. By recalling the result of Theorem~\ref{KdVth1} for HBVMs, in Table~\ref{KdVtab} we list the numerical results obtained by solving the resulting semi-discrete problem (\ref{KdV1}) with time-step $h=10/n$, in terms of: execution time (in sec); maximum solution and Hamiltonian errors, $e_u$ and $e_H$, respectively;  rate of convergence, where appropriate.\footnote{Also in this case, for the Gauss method a super-convergence occurs in the Hamiltonian error.} For the SHBVM method, we also list the used values of $k$ and $s$, the latter obtained by using  $tol\approx 10^{-1}\sqrt{\varepsilon}$ in (\ref{tol}) and $k$ suitably larger than $s$. From the obtained results, one sees that the energy-conserving and symplectic methods are almost equivalent, with the higher-order methods performing better than the lower-order ones. Also in this case, however, the spectral method turns out to be the most effective, being able to use much larger time-steps, with uniformly small solution and Hamiltonian errors.

\begin{table}[t]
\caption{Numerical solution of the Korteweg--de Vries equation (\ref{KdVe}) using a time-step $h=10/n$.}\label{KdVtab}
\centering
\begin{tabular}{rccccc}
\toprule
\multicolumn{6}{c}{Gauss 1}\\
\toprule
$n$~ & time & $e_u$ & rate & $e_H$ & rate\\
\midrule
10000 &  5.1 & 1.10e+00 & --- & 9.38e-07 & ---  \\ 
20000 &  8.6 & 2.75e-01 & 2.0 & 5.85e-08 & 4.0  \\ 
30000 & 13.8 & 1.22e-01 & 2.0 & 1.16e-08 & 4.0  \\ 
40000 & 17.6 & 6.88e-02 & 2.0 & 3.67e-09 & 4.0  \\ 
50000 & 20.7 & 4.40e-02 & 2.0 & 1.50e-09 & 4.0  \\ 
\bottomrule
\toprule
\multicolumn{6}{c}{Gauss 2}\\
\toprule
$n$~ & time & $e_u$ & rate & $e_H$ & rate\\
\midrule
10000 &  9.4 & 9.33e-05 & --- & 7.30e-12 & ---  \\ 
20000 & 18.1 & 5.84e-06 & 4.0 & 1.99e-13 & 5.2  \\ 
30000 & 25.3 & 1.15e-06 & 4.0 & 2.84e-13 & ***  \\ 
40000 & 31.3 & 3.65e-07 & 4.0 & 6.54e-13 & ***  \\ 
50000 & 39.5 & 1.50e-07 & 4.0 & 7.25e-13 & ***  \\ 
\bottomrule
\toprule
\multicolumn{6}{c}{HBVM(2,1)}\\
\toprule
$n$~ & time & $e_u$ & rate & $e_H$ \\
\midrule
10000 &  9.4 & 9.65e-01 & --- & 6.39e-14   \\ 
20000 & 16.6 & 2.42e-01 & 2.0 & 5.68e-14   \\ 
30000 & 21.8 & 1.08e-01 & 2.0 & 7.11e-14   \\ 
40000 & 29.0 & 6.05e-02 & 2.0 & 6.39e-14   \\ 
50000 & 33.6 & 3.87e-02 & 2.0 & 6.39e-14   \\ 
\bottomrule
\toprule
\multicolumn{6}{c}{HBVM(3,2)}\\
\toprule
$n$~ & time & $e_u$ & rate & $e_H$ \\
\midrule 
10000 & 12.8 & 8.49e-05 & --- & 5.68e-14   \\ 
20000 & 24.4 & 5.32e-06 & 4.0 & 5.68e-14   \\ 
30000 & 34.6 & 1.05e-06 & 4.0 & 6.39e-14   \\ 
40000 & 43.0 & 3.32e-07 & 4.0 & 5.68e-14   \\ 
50000 & 54.5 & 1.36e-07 & 4.0 & 7.11e-14   \\ 
\bottomrule
\multicolumn{6}{c}{SHBVM}\\
\toprule
$n$~ & time & $k$ & $s$ & $e_u$ & $e_H$ \\
\midrule
 400 &    7.2 & 20 & 18 &  1.31e-11 &  4.26e-14 \\
 600 &   7.4  & 16 & 14 &   3.70e-12 & 4.26e-14 \\
 800 &   8.6  & 14 & 12 &   4.75e-12 & 4.26e-14 \\
\bottomrule
\end{tabular}
\end{table}

\subsection{A few remarks} From the obtained results, we can draw a few conclusions, which we report in the sequel.
\begin{description}
\item[Energy-conservation.] When the conservation of energy is not an issue, the performance of energy-conserving HBVMs seems to be comparable with that of the symplectic Gauss formulae of the same order. Clearly, things may change when energy-conservation is an important feature (see, e.g., the example in \cite[Section\,7]{BFCI2015}).

\item[Order of the methods.] From the numerical results, one clearly sees that the second-order methods are outperformed by higher-order HBVMs and/or Gauss methods.  In particular, for problems (\ref{NLSEex}) and (\ref{KdVe}), the second-order HBVM(2,1) method is exactly energy-conserving, and can be regarded as a high-performance implementation of the AVF method in \cite{QMcL2008}. Despite this, its performance is not comparable with that of the higher-order methods.

\item[Spectral methods in time.] The obtained numerical results further confirm what recently observed in \cite{BMR2018,BIMR2018,ABI2018}, i.e., that the use of HBVMs as spectral methods in time is a very promising way of getting very high-performance ODE solvers, due to the effectiveness of the underlying blended iteration described in Section~\ref{HBVMsec}.

\end{description}

\section{Conclusions}\label{endsec}

In this paper we have reviewed the basic facts concerning the use of energy-conserving line integral methods for efficiently solving Hamiltonian PDEs. This has been done by performing, at first, a suitable space discretization, along a Fourier orthonormal basis, thus obtaining a corresponding high-dimensional Hamiltonian problem. In particular, we have studied the semilinear wave equation, the nonlinear Schr\"odinger equation, and the Korteweg--de Vries equation in one dimension. It is worth mentioning, however, that:  as sketched in Section~\ref{hisp}, the used space discretization can be straightforwardly extended to the case of more space dimensions; additional Hamiltonian PDEs have been considered in \cite{BZL2018,BGZ2018}. In the future, we plan to further investigate Hamiltonian PDEs within the same framework.

\vspace{6pt} 

%\reftitle{References}

%%%%%%%%%%%%%%%%%%%%%%%%%%%%%%%%%%%%%%%%%%
\end{document}